\documentclass[10pt]{article}
\usepackage{amsfonts, latexsym,amscd,graphicx,amssymb}

\makeatletter \@addtoreset{equation}{section} \makeatother

\def\text#1{\mbox{\rm #1}}
\newcommand{\be}{\begin{equation}}
\newcommand{\ee}{\end{equation}}
\newcommand{\ba}{\begin{eqnarray}}
\newcommand{\ea}{\end{eqnarray}}

\newcommand{\pa}{\partial}
\newcommand{\f}{\frac}

\begin{document}

\title{ A geometric Birkhoffian formalism for nonlinear RLC networks}
\author{Delia Ionescu,\\
{\small {\it Institute of Mathematics of the Romanian Academy}}\\
{\small {\it P.O. Box 1-764, RO-014700, Bucharest,
 Romania,Delia.Ionescu@imar.ro}}}
 \date{}
 \maketitle
\begin{abstract}
The aim of this paper is to give  a formulation of the dynamics of
nonlinear RLC  circuits as a geometric Birkhoffian system and to
discuss in this context the concepts of regularity,
conservativeness, dissipativeness.  An  RLC circuit, with no
assumptions placed on its topology, will be described by a family
of Birkhoffian  systems, parameterized by a finite number of real
constants which correspond to initial values of certain state
variables of the circuit. The  configuration space and a  special
Pfaffian form, called Birkhoffian, are obtained from the
constitutive relations of the  resistors, inductors and capacitors
involved and from Kirchhoff's laws. Under certain assumptions  on
the voltage-current characteristic for resistors, it is shown that
a Birkhoffian system associated to an RLC circuit is dissipative.
For RLC networks which contain a number of pure capacitor loops or
pure resistor loops the Birkhoffian associated is never regular. A
procedure  to reduce the original configuration space
 to a lower dimensional one, thereby regularizing
the Birkhoffian, it is also presented. In order to illustrate the
results, specific examples are discussed in detail.

\textit{Keywords}: geometric methods in differential equations,
Birkhoffian differential systems, Birkhoffian vector fields,
dissipative dynamical systems, electrical networks

MSC: 34A26, 58A20, 94C

\end{abstract}

\section{Introduction}

Lagrangian and Hamiltonian mechanics continue to attract a large
amount of attention in the literature, because many mechanical and
electro-mechanical systems may be modelled within these
frameworks.  During the past ten years, a far-reaching
generalization of the Hamiltonian framework has been developed in
a series of papers. This generalization, which is based on the
geometric notion of generalized Dirac structure, gives rise to
implicit Hamiltonian systems (see for example papers by Maschke
and  van der Schaft \cite{maschke2}, van der Schaft
\cite{schaft}).
 In the papers by Yoshimura and  Marsden
 \cite{yoshimura}
the concept of Dirac structures and variational principle are used
to define and develop the basic properties of implicit Lagrangian
systems.

An alternative approach to the study of dynamical systems which
appears to cover a wide class of systems, among them the
nonholonomic systems, the degenerate systems and the dissipative
ones, is the Birkhoffian formalism, a global formalism of the
dynamics of implicit systems of second order ordinary differential
equations on a manifold. The classical book by G. D. Birkhoff
\cite{bir}, contains in Chapter I many interesting ideas about
classical dynamics from the viewpoint of differential geometry. In
order to present  these ideas in a coordinate free fashion, one
considers   the formalism of 2-jets (see for example  Kobayashi
and Oliva \cite{oliva}). The space of configurations is a smooth
$m$-dimensional differentiable connected manifold  and the
covariant character of the Birkhoff generalized forces is obtained
by introducing the notion of elementary work, called Birkhoffian,
a special Pfaffian form defined  on the 2-jets manifold.   The
dynamical system associated to this Pfaffian form is a subset of
the 2-jets manifold which defines an implicit second order
ordinary differential system. The notion of Birkhoffian allows to
formulate the  concepts of reciprocity, regularity, affine
structure in the accelerations, conservativeness, in an intrinsic
way.

The   electrical circuits theory benefits from many tools
developed in mathematics. In order to study  the dynamics of LC
and RLC electrical circuits,
various Lagrangian and Hamiltonian formulations have been considered in the 
literature (see for example  \cite{bloch}, \cite{moser},
\cite{maschke}, \cite{maschke2},  \cite{schaft}, \cite{smale}, and
the references therein).
 To
describe LC circuit equations, Hamiltonian formulations have been
used more often. In \cite{maschke}, the dynamics of a nonlinear LC
circuit
is shown to be of Hamiltonian nature with respect to a certain  Poisson 
bracket which may be degenerate, that is, non-symplectic.
 The dynamics of  ''complete''   RLC networks was
described by Brayton and Moser \cite{moser} in terms of a function
of  inductor currents and capacitor voltages, called the mixed
potential function.
 However, for all those formulations,  a certain topological assumption on the electrical circuit
appears to be crucial, that is, the 
circuit is supposed to  contain neither loops of capacitors
 nor cutsets of inductors.
In  \cite{maschke2}, \cite{schaft} and \cite{bloch} the Poisson bracket is replaced by
the more general notion of a Dirac structure 
on a vector space, leading to implicit Hamiltonian systems. In
this formalism, it is possible to include LC networks which do not
obey the topological assumption mentioned before. In
\cite{yoshimura} an example of  LC circuit  in the context of
implicit Lagrangian systems is given  for a degenerate Lagrangian
system with holonomic constraints.

\noindent The potential relevance of the Birkhoffian formalism in
the context of electrical circuits is discussed  by Ionescu and
Scheurle \cite{io}, where  a formulation of general nonlinear LC
circuits within the framework of Birkhoffian dynamical systems
on manifolds is presented. In \cite{io}  specific 
examples of  electrical networks are discussed in this framework.
These are networks which contain closed loops formed by capacitors, as well as 
inductor cutsets, and also LC networks which contain independent
voltage sources as well as independent current sources.

In the paper at hand we present a formulation of the dynamics of
nonlinear RLC electrical circuits  within the framework of
Birkhoffian systems. Based on   Kirchhoff's laws and the
constitutive relations for the  resistors, inductors, capacitors
involved, we get for a nonlinear RLC circuit, a whole family of
configuration spaces and a  special Pfaffian forms, called
Birkhoffians. The configurations spaces are parameterized by a
finite number of real constants which correspond to initial values
of certain state variables of the circuit.  In particular, we can
allow pure capacitor loops as well as pure inductor cutsets. Under
certain assumptions  on the voltage-current characteristic for
resistors, it is shown that a Birkhoffian system associated to an
RLC circuit is dissipative.

The  paper is organized as follows. In Section 2, we recall the basics of Birkhoffian systems (see \cite{bir}) 
presented from the viewpoint of differential geometry using the
formalism of jets (see \cite{oliva}).
In particular, we introduce the notion of  dissipative Birkhoffian system, in order to be able to treat the case of RLC networks 
 in the next section.
In Section 3, our Birkhoffian formulation of the dynamic equations
of a nonlinear RLC circuit is given. Properties of the
corresponding Birkhoffian such as its regularity and its
dissipativness are also discussed in this section. For electrical
RLC networks which contain a number of only capacitor loops or only resistors  loops, we present a systematic 
procedure to reduce the original configuration  space  to a lower dimensional one, thereby regularizing the 
Birkhoffian. On the reduced configuration space the reduced
Birkhoffian  will still be dissipative.
 Finally, in section 4 we consider two specific examples.  These examples are supposed to serve our purpose of demonstrating 
the power of the Birkhoffian approach.

\textit{\textbf {Acknowledgement.}} I would like to express my
gratitude  to Professor Marsden for  valuable discussions
concerning the topic of this paper.

\section{Birkhoffian systems}

For a smooth m-dimensional differentiable connected manifold $M$,
we consider the tangent bundles $(TM,\pi_M,M)$ and ($TTM, 
\pi_{TM}, TM$).  Let $q=(q^1$, $q^2,$..., $q^m)$ be a local
coordinate system on $M$. This induces   natural local coordinate
systems  on $TM$ and  $TTM$, denoted by ($q,\dot{q}$),
respectively ($q,\, \dot{q},\,dq,\, d\dot{q}$). The \textit{2-jets
manifold} $J^2(M)$ is a
 $3m$-dimensional submanifold of $TTM$ defined by
 \be
J^2(M)=\left\{ z\in TTM \,\, | \,\,
T\pi_{M}(z)=\pi_{TM}(z)\right\} \label{0} \ee where
$T\pi_{M}:TTM\to TM$ is the tangent map of $\pi_{M}$. We write
$\pi_J:=\pi_{TM}\arrowvert _{J^2(M)}=T\pi_M\arrowvert_{J^2(M)}$.
($J^2(M),\, \pi_J,\, TM$), called the \textit{2-jet bundle} (see \cite{oliva}), is an affine bundle modelled on the vertical vector 
bundle ($V(M),\, \pi_{TM}\arrowvert_{V(M)},\, TM$),  $
V(M)=\bigcup_{v\in TM} V_v(M) $ , where
$V_v(M)=\{z\in T_vTM\, \arrowvert\, (T\pi_M)_v(z)=0\}$. In \cite{marsden}, \cite{saunders} this bundle is denoted by $T^2(M)$ and 
named \textit{second-order tangent bundle}. In natural local
coordinates, the equality in (\ref{0}) yields
 $(q,\, \dot{q},\,\dot{q},\, d\dot{q}\arrowvert_{J^2(M)})$ as a local coordinate system on $J^2(M)$. We set $\ddot{q}:=d\dot 
{q}\arrowvert_{J^2(M)}$.
 Thus, a local coordinate system $q$ on $M$ induces the  natural local  coordinate system
 $(q,\, \dot{q},\,\ddot{q})$ on $J^2(M)$. For further details on this affine bundle see  \cite{marsden}, \cite{oliva}, 
\cite{saunders}.

A \textit{\textbf {Birkhoffian}} corresponding to the configuration manifold $M$ is a smooth 1-form $\omega$ on $J^2(M)$ such that, 
for any $x\in M$, we have \be i_{x}^*\omega=0\label{defbir} \ee
where $i_x:\beta^{-1}(x)\to J^2(M)$ is the embedding of the submanifold $\beta^{-1}(x)$ into $J^2(M)$, $\beta=\pi_M\circ \pi_J$. 
From this definition it follows that, in the  natural local coordinate system ($q,\, \dot{q},\, \ddot{q}$) of $J^2(M)$, a 
Birkhoffian $\omega$ is given by \be \omega=\sum^{m}_{j=1}Q_j(q,\,
\dot{q},\, \ddot{q})dq^j\label{bircor} \ee with certain functions
$Q_j:J^2(M)\to \mathbf{R}$.

 \noindent The pair ($M,\, \omega$) is said to be a \textit{ \textbf
{Birkhoff system}} (see \cite{oliva}).

\noindent The \textit{\textbf{differential system associated to a
Birkhoffian}} $\omega $ (see \cite{oliva} ) is  the set (maybe
empty)  $D(\omega$), given by \be D(\omega):=\left\{z\in J^2(M)\
\arrowvert\, \omega(z)=0 \right\}\label{difsistem'} \ee The
manifold $M$ is the \textit{space of configurations} of
$D(\omega)$, and $D(\omega)$ is said to have $m$ 'degrees of
freedom'. The $Q_i$ are the 'generalized external forces'
associated to the local coordinate system $(q)$. In the natural
local coordinate system, $D(\omega)$ is characterized by the
following implicit system of second order ODE's \be Q_j(q,\,
\dot{q},\, \ddot{q})=0 \textrm{ for all }
j=\overline{1,m}\label{difsistem} \ee \noindent {\it We conclude
that the Birkhoffian formalism is  a global formalism for the
dynamics of implicit systems of second order  differential
equations on a manifold}.

 \noindent Let us now associate a vector
field to a Birkhoffian $\omega$.\\
 A \textit{vector field} $Y$ on the manifold $TM$ is a smooth
function $Y:TM\to TTM$ such that $\pi_{TM}\circ Y$=id.
Any vector field $Y$ on $TM$ is called  \textit{second order vector field on TM} if and only if $T{\pi_M}(Y_v)=v$ for all $v\in 
TM$.\\
 A cross section $X$ of the affine bundle ($J^2(M),\, \pi_J,\, TM$), that is, a smooth function $X:TM\to J^2(M)$ such that  
$\pi_J\circ X$=id,
 can be identified with a special vector field on $TM$, namely, the second order vector field on $TM$ associated to $X$.  Indeed, because ($J^2(M),\, 
\pi_J,\, TM$) is a sub-bundle of
($TTM, \pi_{TM}, TM$) as well as of ($TTM, T\pi_{M}, TM$), its sections can be regarded as sections of these two tangent bundles. Thus, using 
the canonical   embedding  $i:J^2(M)\to TTM$,  $X$ can be
identified with $Y$, that is,  $Y=i\circ X$.\\
In  natural  local coordinates, a second order vector field can be
represented as \be Y=\sum^{m}_{j=1}\left[\dot{q}^j\f{\pa}{\pa
q^j}+\ddot{q}^j(q,\, \dot{q})\f{\pa}{\pa
\dot{q}^j}\right]\label{vectorfield} \ee

A \textit{\textbf {Birkhoffian vector field}} associated to a Birkhoffian $\omega$ of $M$ (see \cite{oliva}) is a smooth second 
order vector field on $TM$, $Y=i\circ X$, with $X:TM\to J^2(M)$,
such that $ImX\subset\, D(\omega)$, that is \be X^*\omega=0 \ee
 In the natural  local coordinate system, a Birkhoffian vector field is given by the expression (\ref{vectorfield}), such 
that $Q_j(q,\dot{q},\ddot{q}(q,\dot{q}))=0$.

A Birkhoffian $\omega
$ is \textit{\textbf {regular}} if and only if \be
\textrm{det}\left[\f{\pa Q_j}{\pa \ddot{q}^i}(q,\, \dot{q},\,
\ddot{q})\right]_{i,j=1,...,m}\neq 0\label{regular} \ee
for all $(q,\, \dot{q},\, \ddot{q})$, and for each $(q,\, \dot{q})$, there exists $\ddot{q}$ such that $Q_j(q,\, \dot{q},\, 
\ddot{q})=0,\, j=1,...,m.$

\noindent If a Birkhoffian $\omega$ of $M$ is regular, then it satisfies \textit{the principle of determinism}, that is, there exists an 
unique Birkhoffian vector field $Y=i\circ X$ associated to
$\omega$ such that $Im\, X=D(\omega)$ (see \cite{oliva}).

 A Birkhoffian $\omega$ of the configuration space $M$ is called \textit{\textbf{conservative}} if and only if there exists a smooth function 
$E_{\omega}:TM\to \mathbf{R}$ such that

\be
 (X^*\omega)Y=dE_{\omega}(Y)\label{conserv'}
\ee for all second order vector fields $Y=i\circ X$ (see
\cite{oliva}).

\noindent Equation (\ref{conserv'}) is  equivalent, in the natural
local coordinate system, to the identity (see \cite{bir}, p. 16,
eq. 4) \be
\sum^{m}_{j=1}Q_j(q\, \, \dot{q},\, \ddot{q})\dot{q}^j=\sum^{m}_{j=1}\left[\f{\pa E_{\omega}}{\pa q^j}\dot{q}^j+\f{\pa E_{\omega}}{\pa 
\dot{q}^j}\ddot{q}^j\right]\label{conserv} \ee $E_{\omega}$ is
constant on $TM$ if and only if $dE_{\omega}(Y)=0$ for all second
order vector fields $Y$ on $TM$ (see \cite{oliva}).\\
 If $\omega$ is conservative and $Y$ is a Birkhoffian
vector field, then (\ref{conserv'}) becomes \be dE_{\omega}(Y)=0
\ee This means that $E_{\omega}$ is constant along the
trajectories of $Y$.

 We  now introduce the concept of a dissipative
Birkhoffian.\\
 A \textit{vertical 1-form on $TM$} (see for example
\cite{yano}) is a 1-form $\Psi$  on $TM$ such that $\Psi(V^v)=0$,
for all $V$ vector field on $M$, where $V^v$ is the vertival lift
of the vector field $V$ to $TM$. The local expression of a
vertical 1-form is \be
\Psi=\sum^{m}_{j=1}\psi_j(q,\dot{q})dq^j\label{vform} \ee
 A  1-form $D$ on $TM$ is called
\textit{\textbf{dissipative}} if and only if $D \textrm{ is
vertical }$ and $D(Y)> 0$, for all $Y$ second order vector field
on $TM$. Allowing for (\ref{vform}), the local expression of $D$
is \be D=\sum^{m}_{j=1}D_{j}(q,\dot{q})dq^j\label{store2} \ee and
from (\ref{store2}), (\ref{vectorfield}), the inequality $D(Y)> 0$
becomes \be \sum^{m}_{j=1}D_{j}(q,\dot{q})\dot{q}^j>0\label{store}
\ee A Birkhoffian $\omega$ of the configuration space $M$ is
called \textit{\textbf{dissipative}}  if and only if there exists
a smooth function $E_{0_{\omega}}:TM\to \mathbf{R}$ such that \be
(X^*\omega)Y = dE_{0_\omega}(Y)+ D(Y) \label{store1} \ee for all
second order vector fields $Y=i\circ X$ on $TM$,  $D$ being a
dissipative 1-form on $TM$.

\noindent Equation (\ref{store1}) is  equivalent, in a local
coordinate system, to the identity \be \sum^{m}_{j=1}Q_j(q\, \,
\dot{q},\, \ddot{q})\dot{q}^j=\sum^{m}_{j=1}\left[\f{\pa
E_{0_\omega}}{\pa q^j}\dot{q}^j+\f{\pa E_{0_\omega}}{\pa
\dot{q}^j}\ddot{q}^j+ D_{j}(q,\dot{q})\dot{q}^j
\right]\label{dissip} \ee

\noindent In view of (\ref{store}), we obtain from (\ref{store1}),
\be
 (X^*\omega)Y> dE_{0_{\omega}}(Y)\label{storage}
\ee for all second order vector field $Y=i\circ X$. That is
equivalent, in local coordinates, to the dissipation inequality
\be \sum^{m}_{j=1}Q_j(q\, \, \dot{q},\, \ddot{q})\dot{q}^j>
\sum^{m}_{j=1}\left[
 \f{\pa E_{0{\omega}}}{\pa q^j}\dot{q}^j+\f{\pa E_{0_{\omega}}}{\pa
 \dot{q}^j}\ddot{q}^j\right]
.\label{storage'}
\ee
 If $\omega$ is a dissipative Birkhoffian  and $Y$ is the
Birkhoffian vector field, then (\ref{storage}) becomes \be
dE_{0_\omega}(Y)< 0. \ee This means that $E_{0_\omega}$ is
nonincreasing along the trajectories of $Y$.

If the dissipative 1-form on $TM$ has the particular expression
\be D=\sum^{m}_{i,j=1}\mathcal{D}_{ij}(q)\dot{q}^jdq^i \ee when we
calculate the function  $D(Y)$ on $TM$,  we obtain the so called
\textit{Rayleigh dissipation function} $\mathcal{R}:TM \to
\mathbf{R}$ \be
\mathcal{R}(q,\dot{q})=\sum^{m}_{i,j=1}\mathcal{D}_{ij}(q)\dot{q}^i\dot{q}^j
\ee

\section{RLC circuit dynamics}

A simple electrical circuit provides us with an \textit{oriented
connected} graph, that is, a collection of points, called nodes,
and a set of connecting lines or arcs, called branches, such that
in each branch is given a direction and there is at least one path
 between any two
nodes. A path is a sequence of branches such that the origin of
the next branch coincides with the end of the previous one. The
graph will be assumed to be \textit{planar}, that is, it can be
drawn in a plane without
branches crossing. For the graph theoretic terminology, see, for 
example \cite{graph}.\\
 Let $b$  be the total number of branches in the graph,
$n$ be one less than the number of  nodes and $m$
be the cardinality of a 
selection of loops that cover the whole graph. Here, a loop is a
path such that the first and last node coincide and that does not
use the same branch  more than once. By Euler's polyhedron
formula, $b=m+n$.

\noindent
 We choose a reference node and a  current direction in
each   $l$-branch of the graph,
 $l=1,...,b$. We also consider a covering of the graph with $m$ loops,
   and a current direction in  each $j$-loop, 
$j=1,...,m$. We assume that the associated graph has at least one
loop, meaning that $m>0$.\\
 A graph  can be described by matrices: a ($bn$)-matrix
$B\in \mathfrak{M}_{bn}(\mathbf{R})$, rank$(B)=n$, called
\textit{incidence matrix} and a ($bm$)-matrix $A\in \mathfrak{M}_
{bm} (\mathbf{R})$, rank$(A)=m$, called \textit{loop matrix}.
These matrices contain only 0,1, -1. An element of the matrix $B$
is 0 if a branch  $b$ is not incident with a node $n$, 1 if branch
$b$ enters  node  $n$ and -1 if
 branch $b$ leaves node $n$, respectively. An element of the matrix $A$ is 0 if a branch
 $b$ does
not belong to a loop $m$, 1 if  branch $b$
 belongs to loop $m$
and their directions agree
 and -1 if
 branch $b$ belongs to loop $m$  and their  directions
oppose, respectively. For the fundamentals of electrical circuit
theory, see, for example \cite{chua}.

\noindent The states of the circuit have two components, the
currents through the branches, denoted by $\textsc{i}\in
\mathbf{R}^b$,  and the voltages across the branches, denoted by
$v\in \mathbf{R}^b$.
 Using the matrices
$A$ and $B$, Kirchhoff's current law and Kirchhoff's voltage law
can be expressed by the equations \be B^T \textsc{i} =0 \quad
(KCL)\label{5} \ee \be A^T v =0\quad (KVL)\label{6} \ee \noindent
Tellegen's theorem establishes a relation between the matrices
$A^T$ and $B^T$: \textit{the kernel of the matrix $B^T$ is
orthogonal to the kernel of the matrix $A^T$}(see e.g.,
\cite{moser} page 5).

 The next step is to introduce the branch elements in
a simple electrical circuit. The branches of the
 graph associated to a  RLC electrical circuit, can be classified
into three categories: resistive branches,  inductor branches,
 and capacitor branches. Let $r$ denote the number of resistive branches,
 $k$ the number of inductor branches and $p$ the number of capacitor
 branches, respectively. We assume that just one electrical device is associated to each
branch, then, we have $b=r+k+p$. Thus, we can write
$\textsc{i}=(\textsc{i}_{[\Gamma]},
\textsc{i}_{(a)},\textsc{i}_{\alpha})\in \mathbf{R}^r\times
\mathbf{R}^k\times \mathbf{R}^p\simeq\mathbf{R}^b$, where
$\textsc{i}_{[\Gamma]}$, $\textsc{i}_ {(a)}$,
$\textsc{i}_{\alpha}$ are the currents through the resistors, the
inductors,  the capacitors, respectively, and
$v=(v_{[\Gamma]},v_{(a)},v_{\alpha})\in \mathbf{R}^r\times
\mathbf{R}^k\times\mathbf{R}^p\simeq\mathbf{R}^b$, where
$v_{[\Gamma]}$, $v_{(a)}$, $v_\alpha$ describe the voltage drops
across the resistors,  the inductors, the capacitors,
respectively.\\
Each capacitor is supposed to be charge-controlled. For the
nonlinear capacitors we assume \be
v_\alpha=C_\alpha(\textsc{q}_{\alpha}), \quad \alpha=1,..., p
\label{2} \ee where the functions
$C_\alpha:\mathbf{R}\longrightarrow \mathbf{R}\backslash \{0\}$
are smooth and invertible, and the 
$\textsc{q}_\alpha$'s denote the charges of the capacitors.
 The current through a  capacitor is given by  the
time-derivative
 of the corresponding  charge
\be \textsc{i}_\alpha=\f{d\textsc{q}_{\alpha}}{dt}, \quad
\alpha=1,...,p \label{4'} \ee $t$ being the time variable.\\
 Each inductor is supposed to be current-controlled. For
the  nonlinear inductors we assume \be
v_a=L_{a}({\textsc{i}_a})\f{d{\textsc{i}_a}}{dt}, \quad a=1,...,k
\label{4} \ee where $L_a:\mathbf{R}\longrightarrow
\mathbf{R}\backslash \{0\}$ are smooth invertible functions.\\
 There are several types of nonlinear resistors, among
them current controlled resistors and  voltage controlled
resistors. Generally, their constitutive relations   are defined
by some continuous functions of $\textsc{i}_{\Gamma}$ and
$v_{\Gamma}$, that is, \\
$
f_{\Gamma}(\textsc{i}_{\Gamma},v_{\Gamma})=0, \Gamma=1,...,r $.

\subsection{Current controlled resistors}

\noindent We first consider the case that the  nonlinear resistors
are current controlled, that is, the constitutive relations are
given by \be v_{\Gamma}=R_{\Gamma}(\textsc{i}_{\Gamma}), \quad
\Gamma=1,...,r\label{2rezistor} \ee where
$R_{\Gamma}:\mathbf{R}\longrightarrow \mathbf{R}$ are smooth
functions. In order to obtain  a dissipative Birkhoffian, we
assume that, for all $x \neq 0$, \be R_{\Gamma}(x)x>0, \quad
\forall \Gamma=1,...,r\label{assumption} \ee that is, for each
nonlinear resistor, the graph of the function $R_{\Gamma}$ lies in
the union of the first and  the third quadrant.\\
For linear resistors  the relations  (\ref{2rezistor}) can be
written in the form \be v_{\Gamma}=\textrm {\scriptsize
R}_{\Gamma}\textsc{i}_{\Gamma}\label{15rez} \ee \noindent where
$\textrm {\scriptsize R}_{\Gamma}> 0$ are real constants.\\
Taking into account (\ref{2}), (\ref{4'}), (\ref{4}),
(\ref{2rezistor}),  the equations (\ref{5}), (\ref{6})  become \be
\left\{
\begin{array}{ll}
B^T
\left(
\begin{array}{c}
\textsc{i}_{\Gamma}\\
\textsc{i}_{a} \\
\f{d\textsc{q}_{\alpha}}{dt}
\end{array}
\right)=0\\
\\
A^T
\left(
\begin{array}{c}
R_{\Gamma}(\textsc{i}_{\Gamma})\\
L_{a}(\textsc{i}_a)\,
\f{d\textsc{i}_a}{dt} \\
C_\alpha(\textsc{q}_{\alpha})
\end{array}
\right)=0
\end{array}
\right.\label{7} \ee \textit{In the following we give a
Birkhoffian formulation for the network described by the system of
equations (\ref{7}), using the same procedure as in \cite{io}.
That is,
 using the first set of equations (\ref{7}),
 we are going to define a  family of
  $m$-dimensional affine-linear configuration spaces
  $M_c\subset\mathbf{R}^b$ parameterized by a constant vector
  $c$ in $\mathbf{R}^n$. This vector is related to the initial
  values of the $\textsc{q}$-variables at some instant of time.
A Birkhoffian $\omega_c$ on the configuration space $M_c$ arises
from a linear combination of the second set of equations
(\ref{7}). Thus, ($M_c,\omega_c$) will be a family of Birkhoff
systems that describe the considered RLC network.}

\noindent We  notice that the first set of equations (\ref{7})
 remains exactly the same for linear and nonlinear electrical devices.
 Thus, for obtaining the  configuration space, it is not important
 whether the devices are linear or nonlinear.

\noindent Let  $H: \mathbf{R}^b\longrightarrow\mathbf{R}^n$
 be the linear map  that, with respect to a  coordinate system
 ($x^1,..., x^b$) on
$\mathbf{R}^b$, is given  by

\be H(x^1,...,x^b)=B^T \left(
\begin{array}{c}
x^1 \\
\vdots\\
x^b
\end{array}
\right)\label{9} \ee Then, $H^{-1}(c)$, with $c$ a constant vector
in $\mathbf{R}^n$, is an affine-linear  subspace in
$\mathbf{R}^b$. Its dimension is $m=b-n$, because rank$(B)=n$.

 \textit{We
define  $M_c$ as }\be M_c:=H^{-1}(c)\label{8} \ee We denote local
coordinates  on $M_c$ by $q=(q^1,..,q^m)$. Then, the natural
coordinate system on the 2-jets bundle $J^2(M_c)$ is given by
$(q,\dot{q},\ddot{q})$. \\
We will now represent the  Birkhoffian in a specific coordinate
system
on $M_c$:\\
In the vector spaces $\mathbf{R}^r$, $\mathbf{R}^k$, we identify
points and vectors \be
\textsc{i}_{\Gamma}:=\f{d\textsc{q}_{[\Gamma]}}{dt},\quad
\textsc{i}_{a}:=\f{d\textsc{q}_{(a)}}{dt}\label{var} \ee with
$(\textsc{q}_{[\Gamma]})_{\Gamma=1,...,r}$,
 $(\textsc{q}_{(a)})_{a=1,...,k}$
 coordinate systems on $\mathbf{R}^r$, respectively, on $\mathbf{R}^k$.
Taking into account (\ref{var}) and the fact that the matrix $B^T$
is a constant matrix, we integrate  the first set of equations
(\ref{7}) to arrive at  \be B^T \left(
\begin{array}{c}
\textsc{q}_{[\Gamma]} \\
\textsc{q}_{(a)} \\
\textsc{q}_{\alpha}
\end{array}
\right)=c\label{7''} \ee with $c$ a constant vector in
$\mathbf{R}^n$.\\
Likewise  consider coordinates  in $\mathbf{R}^b\simeq
\mathbf{R}^r\times\mathbf{R}^k\times\mathbf{R}^p$
\vspace{-0.5cm}\ba
&&x^1:=\textsc{q}_{[1]},\,...,\,x^r:=\textsc{q}_{[r]},\,x^{r+1}:=\textsc{q}_{(1)},\,...
,\,x^{r+k}:=\textsc{q}_{(k)}, \nonumber\\
&&x^{r+k+1}:=\textsc{q}_{1}, \,...,\,x^b:=\textsc{q}_{p}\label{cs}
\ea From (\ref{9}), (\ref{8}), we see that we can define
coordinates on $M_c$ by solving the equations (\ref{7''}) in terms
of an appropriate set of $m$ of the $\textsc{q}$-variables, say
$q=(q^1,...,q^m)$. In other words, we express any of the
$x$-variables as a function of $q=(q^1,...,q^m)$, namely,
\vspace{-0.5cm}\ba
&&x^{\Gamma}=\sum_{j=1}^{m}\mathcal{N}^{\Gamma}_jq^j+const,\,
 \Gamma=1,...,r\nonumber\\
&&x^a=
\sum_{j=1}^{m}\mathcal{N}^{a}_jq^j+const,\,
a=r+1,...,r+k\nonumber\\
&& x^{\alpha}= \sum_{j=1}^{m}\mathcal{N}^{\alpha}_jq^j+const, \,
\alpha=r+k+1,...,b \label{10'} \ea with certain constants
$\mathcal{N}^{\Gamma}_j$, $\mathcal{N}^{a}_j$,
$\mathcal{N}^{\alpha}_j$. Here we can think of the constants
$const$ as being initial values of the $x$-variables at some
instant
of time.\\
From  (\ref{4'}), (\ref{var}),
 (\ref{cs}) and  differentiating  (\ref{10'}) we get
\be \textsc{i}=\mathcal{N}\dot{q}\label{N} \ee with the matrix of
constants $\mathcal{N}\in\mathfrak{M}_{bm} (\mathbf{R})$, for some
$\dot{q}\in \mathbf{R}^m$.\\
Using  Tellegen's theorem and a fundamental theorem of linear
algebra, we now find a relation between the matrices $\mathcal{N}$
and $A$.\\
 By a fundamental theorem of linear algebra we have \be
(Ker(A^T))^\perp=Im(A) \label{linearalgebra}\ee where  $A\in
\mathfrak{M}_{bm}(\mathbf{R})$, $Ker(A^T):=\{x\in \mathbf{R}^b\,
|\, A^Tx=0\}$ is the kernel of $A^T$, $Im(A):=\{x\in
\mathbf{R}^b\, |\, Ay=x,\, \textrm{for some} \, y\in \mathbf{R}^m
\}$ is the image of $A$ and  $^\perp$ denotes the orthogonal
complement in $\mathbf{R}^b$ of
the respective vector subspace.\\
 For the incidence  matrix $B\in
\mathfrak{M}_{bn}(\mathbf{R})$ and the loop matrix $A\in
\mathfrak{M}_{bm} (\mathbf{R})$, which satisfy Kirchoff's law
(\ref{5}), (\ref{6}), Tellegen's theorem writes as \be
Ker(B^T)=(Ker(A^T))^\perp \label{tellegen}\ee From the first set
of equations in (\ref{7}), and by constraction of the matrix
$\mathcal{N}$ in (\ref{N}), we have \be
Ker(B^T)=Im(\mathcal{N})\label{N1}\ee Therefore, using
(\ref{linearalgebra}), (\ref{tellegen}), (\ref{N1}), we obtain
$Im(A)=Im(\mathcal{N})$. Then, another application of
(\ref{linearalgebra}) yields \textit{ \be
Ker(A^T)=Ker(\mathcal{N}^T)\label{ker}\ee} Taking into account
(\ref{ker}), we see that there exists a nonsingular matrix
$\mathfrak{C}\in \mathfrak{M}_{mm}(\mathbf{R})$ satisfing \be
\mathfrak{C}A^T=\mathcal{N}^T\label{equivalence} \ee The matrix
$\mathfrak{C}$ provides a relation between the  vector of the $m$
independent loop currents and the coordinate vector $q$ introduced
on  $M_c$.

 \textit{Taking into account (\ref{ker}), we
define the Birkhoffian $\omega_c$ of $M_c $ such that the
differential system (\ref{difsistem}) is the linear combination of
the second set of equations in (\ref{7}) obtained by replacing
$A^T$ with the matrix $\mathcal{N}^T$. Thus, in terms of
q-coordinates as chosen before, the expressions of the components
$Q_j(q, \, \dot{q},\,\ddot{q})$  are} \be Q_j(q, \,
\dot{q},\,\ddot{q})=F_j(\dot{q})\ddot{q}+ H_j(\dot{q})+G_j(q),
\quad j=1,...,m\label{bir} \ee\textit{where} {\footnotesize \be
\hspace{0cm}F_j(\dot{q})\ddot{q}:=\sum^{r+k}_{a=r+1}\mathcal{N}^a_{j}L_{a-r}\left(\sum^{m}_{l=1}\mathcal{N}^
{a}_{l}\dot{q}^l\right)\left(\sum ^{m}_{i=1}\mathcal{N}^
{a}_{i}\ddot{q}^i\right)=\sum^{m}_{i=1}
\sum^{r+k}_{a=r+1}\mathcal{N}^a_{j}\mathcal{N}^
{a}_{i}\widetilde{L}_{a-r}\left(\dot{q}\right)\ddot{q}^i
\label{bir2}\ee} {\footnotesize \be H_j(\dot{q}):=\sum
^{r}_{\Gamma=1}\mathcal{N}_j^{\Gamma}R_{\Gamma}\left(\sum^{m}_{l=1}\mathcal{N}^
{\Gamma}_{l}\dot{q}^l\right)= \sum
^{r}_{\Gamma=1}\mathcal{N}_j^{\Gamma}\widetilde{R}_{\Gamma}\left(\dot{q}\right)
\label{bir3} \ee} {\footnotesize \be \hspace{0cm}G_j(q):=\sum
^{b}_{\alpha=r+k+1}\mathcal{N}_j^{\alpha}C_{\alpha-r-k}\left(\sum^{m}_{l=1}\mathcal{N}^
{\alpha}_{l}q^l+const\right)=\sum
^{b}_{\alpha=r+k+1}\mathcal{N}_j^{\alpha}\widetilde{C}_{\alpha-r-k}\left(q\right)\label{bir1}
\ee} We note that the Birkhoffian (\ref{bir}) is
\textbf{\textit{not} conservative}. We easily see that there does
not exist a function $E_{\omega}$ such that (\ref{conserv}) is
fulfilled for the Birkhoffian (\ref{bir}), since $\f{\pa^2
E_{\omega}}{\pa q^j\pa \dot{q}^j}\neq \f{\pa^2 E_{\omega}}{\pa
\dot{q}^j\pa q^j}$.

 \textit{For an RLC electrical network with nonlinear
resistors, described by (\ref{2rezistor}), (\ref{assumption}),
nonlinear inductors and capacitors described by (\ref{4}),
respectively (\ref{2}),
 we claim that the Birkhoffian (\ref{bir}) is a
\textbf{ dissipative Birkhoffian}}.

Indeed, in the view of the assumption (\ref{assumption}), the
vertical 1-form $D$ on $TM$ given by \be
D=\sum^{m}_{j=1}H_j(\dot{q})dq^j\label{linearvert} \ee with
$H_j(\dot{q})$ in (\ref{bir3}), is dissipative, that is, \be
\sum_{j=1}^{m}H_j(\dot{q})\dot{q}^j=\sum ^{r}_{\Gamma=1} \left[
R_{\Gamma}\left(\sum_{j=1}^{m}\mathcal{N}_{j}^{\Gamma}\dot{q}^j\right)\right]\left(
\sum_{j=1}^{m}\mathcal{N}_{j}^{\Gamma}\dot{q}^j\right)>0\label{}
\ee
 We showed in \cite{io} that the following  smooth function $E_{0_\omega}$ on $TM$
{\footnotesize \ba \hspace{-0.5cm} E_{0_\omega}=&&
\sum_{a=r+1}^{r+k}\sum_{l=1}^m\sum_{i_1<...<i_l=1}^{m}(-1)^{l+1}\underbrace{\int_{...}
\int}_l \left[ \widetilde{L}_{a-r}^{(l-1)}(\dot{q})\mathcal{N}^{
a}_{i}
\dot{q}^i+\right.\nonumber\\
&&\left.(l-1)\widetilde{L}_a^{(l-2)}(\dot{q})\right]\mathcal{N}^{
a}_{i_1}...\mathcal{N}^{ a}_{i_l}
d\dot{q}^{i_1}...d\dot{q}^{i_l}+\nonumber\\
&&\sum_{\alpha=r+k+1}^{b}\sum_{l=1}^m\sum_{i_1<...<i_l=1}^{m}(-1)^{l+1}\underbrace{\int_{...}
\int}_l
\widetilde{C}_{\alpha-r-k}^{(l-1)}(q)\mathcal{N}^{\alpha}_{i_1}...\mathcal{N}^{\alpha}_{i_l}
dq^{i_1}...dq^{i_l}\label{energie} \ea} satisfies the identity \be
\sum^{m}_{j=1}\left[
F_j(\dot{q})\ddot{q}+G_j(q)\right]\dot{q}^j=\sum^{m}_{j=1}\left[\f{\pa
E_{0_{\omega}}}{\pa q^j}\dot{q}^j+\f{\pa E_{0_{\omega}}}{\pa
\dot{q}^j}\ddot{q}^j \right]\label {articoltrecut}\ee According to
(\ref{articoltrecut}), the Birkhoffian (\ref{bir}) satisfies
(\ref{dissip}) with the function $E_{0_\omega}(q,\dot{q})$ given
by (\ref{energie}) and the dissipative 1-form $D$ given by
(\ref{linearvert}). $\quad \blacksquare$

\noindent Let us now discuss the   regularity of the Birkhoffian
given by (\ref{bir}).

 \textit{If there exists in the
network at least one  loop  that contains only capacitors, or
only resistors, or only  resistors and capacitors, then the
Birkhoffian (\ref{bir}) associated to the network is \textbf{never
regular}}.

 In
\cite{io} we have shown that if there exists at least one  loop in
an LC network that contains only capacitors, then the Birkhoffian
associated to the network is never regular. The Birkhoffian
associated to an RLC network  which  contains  at least one loop
formed  only by resistors  or only by resistors and capacitors
 is never regular as
well. The proof is based on the fact that for the $l$-loop which
does not contain any inductor branches, on the column $l$ of the
matrix $A$ we have \be A^a_{l}=0,\quad a=r+1,...,r+k \label{a}\ee
 For the Birkhoffian (\ref{bir}),  the determinant in
 (\ref{regular}) becomes
\be
\textrm{det}\left[\f{\pa Q_j}{\pa \ddot{q}^i}(q,\, \dot{q},\, 
\ddot{q})\right]_{i,j=1,...,m}=\textrm{det}\left[\sum^{r+k}_{a=r+1}\mathcal{N}_j^{\,
a}\mathcal{N}^ {a}_{\,
i}\widetilde{L}_{a-r}\left(\dot{q}\right)\right]_{i,j=1,...,m}\label{det}
\ee From (\ref{equivalence}), we get
$\mathcal{N}_j^a=\sum_{i_1=1}^{m}\mathfrak{C}_j^{i_1}A_{i_1}^a$,
for any $a=r+1,...,r+k$ and taking into account (\ref{a}), we have
\ba \hspace{-0.5cm}\sum^{r+k}_{a=r+1}\mathcal{N}_j^{\,
a}\mathcal{N}^ {a}_{\,
i}\widetilde{L}_{a-r}\left(\dot{q}\right)&=&\sum_{i_1=1 \atop
i_1\neq l}^m\mathfrak{C}^{i_1}_j\mathfrak{C}^{i_1}_i\left[
\sum_{a=r+1}^{r+k}(A^a_{i_1})^2\widetilde{L}_{a-r}\left(\dot{q}\right)\right]
+\nonumber\\
&&\sum_{i_1<j_1 \atop i_1,j_1\neq
l}^m\left(\mathfrak{C}^{i_1}_j\mathfrak{C}^{j_1}_i+\mathfrak{C}^{i_1}_i\mathfrak{C}^{j_1}_j\right)
\left[
\sum_{a=r+1}^{r+k}A^a_{i_1}A^a_{j_1}\widetilde{L}_{a-r}\left(\dot{q}\right)\right]\nonumber\\
&&\label{determinant} \ea Using basic calculus, the determinant of
the matrix with the elements (\ref{determinant}),   is a linear
combination of determinants having at least two linearly dependent
columns. This shows that the determinant in the right hand side of
(\ref{det}) is equal to zero. Thus, the
 Birkhoffian (\ref{bir}) is not regular.
$\quad\blacksquare$

We now discuss  the question, how to proceed in the case than the
Birkhoffian given by (\ref{bir}) is not regular in the sense of
definition (\ref{regular}).

 \textit{If
 there  exists in the network  $m_1<m$  loops which contain only capacitors,
 all the other loops containing at least an inductor, we can
 \textbf{regularize} the Birkhoffian (\ref{bir}) via \textbf{reduction of the
 configuration space}. The reduced configuration space $\bar{M}_{c}$
 of dimension $m-m_1$, is a linear or a nonlinear subspace of $M_c$,
depending on  whether the capacitors  are linear or nonlinear. We
claim that  the Birkhoffian $\bar{\omega}_c$ on the reduced
configuration space $\bar{M}_c$ is  still a \textbf{dissipative
Birkhoffian}. Under certain conditions on the functions $L_a$,
$a=1,...,k,$ which characterize the inductors, the reduced
Birkhoffian $\bar{\omega}_c$ will be a \textbf{regular
Birkhoffian}}.

\noindent Without loss of generality, we can assume that there is
one loop in the network that contains only capacitors and in the
coordinate system we have chosen \be \mathcal{N}_1^{\Gamma}=0, \,
\Gamma=1,..., r,\quad \mathcal{N}_1^{a}=0, \, a=r+1,...,
r+k\label{ngamma}\ee  Thus, the Birkhoffian components
(\ref{bir}), with (\ref{bir2}), (\ref{bir3}), (\ref{bir1}), are
given by,  $j=2,...,m$, \vspace{-0.5cm}{\footnotesize \ba
 \hspace{-0.75cm}
 Q_1(q,\dot{q},\ddot{q})&=&\sum ^{b}_{\alpha=r+k+1}\mathcal{N}_1^{\alpha}\widetilde{C}_{\alpha-r-k}(q)\nonumber\\
\hspace{-0.75cm}Q_j(q,\dot{q},\ddot{q})&=&\sum^{m}_{i=2}
\sum^{r+k}_{a=r+1}\mathcal{N}^a_{j}\mathcal{N}^
{a}_{i}\widetilde{L}_{a-r}\left(\dot{q}\right)\ddot{q}^i+ \sum
^{r}_{\Gamma=1}\mathcal{N}_j^{\Gamma}\widetilde{R}_{\Gamma}\left(\dot{q}\right)+\sum
^{b}_{\alpha=r+k+1}\mathcal{N}_j^{\alpha}\widetilde{C}_{\alpha-r-k}(q)
\nonumber\\
&&\hspace{-0.75cm}\ea} We note that, according to (\ref{ngamma}),
 $\dot{q}^1$ does  not
appear in any function $\widetilde{R}_{\Gamma}(\dot{q})$,
$\widetilde{L}_{a-r}(\dot{q})$
 and  the terms
$\widetilde{L}_{a-r}(\dot{q})\ddot{q}^1$ do not appear in any of
the Birkhoffian components $Q_2(q,\dot{q}, \ddot{q}),..., Q_m(q,
\dot{q}, \ddot{q})$.\\
 We define the $(m-1)$-dimensional nonlinear space
$\bar{M}_{c}\subset M_c$  by \be \bar{M}_c=\{q\in M_c\,\, | \, \,
\sum
^{b}_{\alpha=r+k+1}\mathcal{N}_1^{\alpha}\widetilde{C}_{\alpha-r-k}(q)=0\}
\ee
 By the implicit function theorem, we obtain a
local coordinate system on the reduced configuration space
$\bar{M}_{c}$. Taking $\bar{q}^1:=q^2$,..., $\bar{q}^{m-1}:=q^m$,
the Birkhoffian has the form
$\bar{\omega}_c=\sum^{m-1}_{j=1}\bar{Q}_jd\bar{q}^j$,  \be
\bar{Q}_j(\bar{q}, \dot{\bar{q}},
\ddot{\bar{q}})=\bar{F}_j(\dot{\bar{q}})\ddot{\bar{q}}+
\bar{H}_j(\dot{\bar{q}})+\bar{G}_j(\bar{q}), \quad \textrm{where}
\label{birred} \ee
\be 
\bar{F}_j(\dot{\bar{q}})\ddot{\bar{q}}:= \sum^{m-1}_{i=1}
\sum^{r+k}_{a=r+1}\mathcal{N}^a_{(j+1)}\mathcal{N}^ {a}_{(i+1)}
L_{a-r}\left(\sum^{m-1}_{l=1} \mathcal{N}^a_{(l+1)}
\dot{\bar{q}}^l\right) \ddot{\bar{q}}^i
\ee
\be \bar{H}_j(\dot{\bar{q}}):= \quad \sum
^{r}_{\Gamma=1}\mathcal{N}_{(j+1)}^{\Gamma}
R_{\Gamma}\left(\sum^{m-1}_{l=1} \mathcal{N}^\Gamma_{(l+1)}
\dot{\bar{q}}^l\right)
 \ee
\be \hspace{-0.8cm} \bar{G}_j(\bar{q}):= \sum
^{b}_{\alpha=r+k+1}\mathcal{N}_{(j+1)}^{\alpha} C_{\alpha-r-k}
\left(\mathcal{N}^\alpha_1f(\bar{q}^1,...,\bar{q}^{m-1})+\sum^{m-1}_{l=1}
\mathcal{N}^\alpha_{(l+1)} \bar{q}^l+const\right)
 \label{bir6'} \ee
 $f:U\subset\mathbf{R}^{m-1}\longrightarrow \mathbf{R}$ being the unique function such that
$f(\bar{q}_0)=q^1_0$, $q^1_0\in \mathbf{R}$, and \be \sum
^{b}_{\alpha=r+k+1}\mathcal{N}_{1}^{\alpha} C_{\alpha-r-k}
\left(\mathcal{N}^\alpha_1f(\bar{q}^1,...,\bar{q}^{m-1})+\sum^{m-1}_{l=1}
\mathcal{N}^\alpha_{(l+1)}
\bar{q}^l+const\right)=0\label{implicit} \ee
 for all $\bar{q}=(\bar{q}^1,..., \bar{q}^{m-1})\in U$, with
$U$ a neighborhood of
$\bar{q}_0=(\bar{q}_0^1,...,\bar{q}_0^{m-1})$.\\
 We will now prove that the Birkhoffian (\ref{birred}) is
dissipative.
 In order to do so, we will show that there exists  a
function $\bar{E}_{0_\omega}(\bar{q},\dot{\bar{q}})$ satisfying
\be \sum^{m-1}_{j=1}\bar{Q}_j(\bar{q}\, \, \dot{\bar{q}},\,
\ddot{\bar{q}})\dot{\bar{q}}^j=\sum^{m-1}_{j=1}\left[\f{\pa
\bar{E}_{0_\omega}}{\pa \bar{q}^j}\dot{\bar{q}}^j+\f{\pa
\bar{E}_{0_\omega}}{\pa \dot{\bar{q}}^j}\ddot{\bar{q}}^j+
\bar{D}_{j}(\bar{q},\dot{\bar{q}})\dot{\bar{q}}^j\right]
\label{dissip'} \ee where
$\bar{D}=\sum_{j=1}^{m-1}\bar{D}_{j}(\bar{q},\dot{\bar{q}})d{\bar{q}}^j$
is a dissipative 1-form on $T\bar{M_c}$.

\noindent We consider the following   1-form on $T\bar{M}_c$  \be
\bar{D}=\sum_{j=1}^{m-1}\sum_{\Gamma=1}^{r}\bar{H}_j(\dot{\bar{q}})
d\bar{q}^j \label{Dbar}\ee In the view of the assumption
(\ref{assumption}), the vertical 1-form (\ref{Dbar}) on
$T\bar{M}_c$ is dissipative, that is, \be \sum ^{r}_{\Gamma=1}
\left[
R_{\Gamma}\left(\sum_{j=1}^{m-1}\mathcal{N}_{(j+1)}^{\Gamma}\dot{\bar{q}}^j\right)\right]\left(
\sum_{j=1}^{m-1}\mathcal{N}_{(j+1)}^{\Gamma}\dot{\bar{q}}^j\right)>0
\ee Therefore, (\ref{dissip'}) is fulfilled if
$\bar{E}_{0_\omega}(\bar{q},\dot{\bar{q}})$ can be chosen in such
a way that \be \sum^{m-1}_{j=1}\left[
\bar{F}_j(\dot{\bar{q}})\ddot{\bar{q}}+\bar{G}_j(\bar{q})\right]\dot{\bar{q}}^j
=\sum^{m-1}_{j=1}\left[\f{\pa \bar{E}_{0_\omega}}{\pa
\bar{q}^j}\dot{\bar{q}}^j+\f{\pa \bar{E}_{0_\omega}}{\pa
\dot{\bar{q}}^j}\ddot{\bar{q}}^j\right]\label{dissip''} \ee is
satisfied. Because of the special form of the terms on the left
side of (\ref{dissip''}), we may assume that
$\bar{E}_{0_\omega}(\bar{q},\dot{\bar{q}})$ is a sum of a function
 depending only on $\bar{q}$, and a function  depending only on
$\dot{\bar{q}}$.
From the theory of 
total differentials, a necessary condition for the existence of
such functions is the fulfillment of the following relations \be
\left\{
\begin{array}{ll}
\f {\partial \bar{\mathcal{F}}_j(\dot{\bar{q}})}{\partial
\dot{\bar{q}}^l }-\f {\partial \bar{\mathcal{F}}_l
(\dot{\bar{q}})}{\partial
\dot{\bar{q}}^j}=0\\
\\
\f {\partial \bar{G}_j(\bar{q})}{\partial \bar{q}^l}-\f {\partial
\bar{G}_l(\bar{q})}{\partial \bar{q}^j}=0
\end{array}
\right.\label{bir5} \ee for any $j,l=1,...,m-1$, where \be
\bar{\mathcal{F}}_j(\dot{\bar{q}}):=\sum^{m-1}_{i=1}\sum^{r+k}_{a=r+1}\mathcal{N}^{\,
a}_{(j+1)}\mathcal{N}^ {a}_{(i+1)} L_{a-r} \left(\sum^{m-1}_{l=1}
\mathcal{N}^a_{(l+1)} \dot{\bar{q}}^l\right) \dot{\bar{q}}^i
\label{bir6} \ee  From (\ref{bir6}), we get: \be \hspace{-0.8cm}\f
{\partial \bar{\mathcal{F}}_j(\dot{\bar{q}})}{\partial
\dot{\bar{q}}^l}=
\sum^{r+k}_{a=r+1}\mathcal{N}^a_{(j+1)}\mathcal{N}^{a}_{(l+1)}\widetilde{L}_{a-r}(\dot{\bar{q}})+
\sum^{m-1}_{i=1}\sum^{r+k}_{a=r+1}\mathcal{N}^a_{(j+1)}\mathcal{N}^
{a}_{(i+1)}\mathcal{N}^
{a}_{(l+1)}\widetilde{L}_{a-r}'(\dot{\bar{q}})\dot{\bar{q}}^i\label{bir7}
\ee where  $\widetilde{L}'_{a-r}:=\f
{d\widetilde{L}_{a-r}(\eta)}{d\eta}$. Then,  the left side of the
first relation in (\ref{bir5}) becomes \vspace{-0.5cm}\ba
\sum^{r+k}_{a=r+1} &&\left[
\left(\mathcal{N}^{a}_{(j+1)}\mathcal{N}^a_{(l+1)}-\mathcal{N}^{a}_{(l+1)}\mathcal{N}^{a}_{(j+1)}\right)\widetilde{L}_{a-r}(\dot{\hat{q}})
+\nonumber\right.\\
&&\left.\left(\mathcal{N}^{a}_{(j+1)}\mathcal{N}^a_{(l+1)}-\mathcal{N}^{a}_{(l+1)}\mathcal{N}^{a}_{(j+1)}\right)
(\sum^{m-1}_{\textrm{\scriptsize
i}=1}\mathcal{N}^a_{(i+1)}\dot{\hat{q}}^{i})\widetilde{L}'_{a-r}(\dot{\hat{q}})\right]
\label{bir10} \ea We  easily see that the expression in
(\ref{bir10}) is zero, thus the first relation in (\ref{bir5}) is
fulfilled.\\
 From (\ref{bir6'}), the second relation in (\ref{bir5})
reads as \ba &&\sum
^{b}_{\alpha=r+k+1}\left\{\mathcal{N}_{(j+1)}^{\alpha}
\widetilde{C}'_{\alpha-r-k}(\bar{q})\left[\mathcal{N}^\alpha_1\f{\pa
f(\bar{q})}{\pa \bar{q}^l}
+\mathcal{N}_{(l+1)}^{\alpha}\right]-\right.\nonumber\\
&& \left. \quad \quad \quad \quad\mathcal{N}_{(l+1)}^{\alpha}
\widetilde{C}'_{\alpha-r-k}(\bar{q})\left[\mathcal{N}^\alpha_1\f{\pa
f(\bar{q})}{\pa \bar{q}^j}
+\mathcal{N}_{(j+1)}^{\alpha}\right]\right\}=0 \label{c}\ea where
$\widetilde{C}'_{\alpha-r-k}:=\f
{d\widetilde{C}_{\alpha-r-k}(\eta)}{d\eta}$. The relation
(\ref{c}) reduces to \be \hspace{-0.7cm}\sum ^{b}_{\alpha=r+k+1}
\mathcal{N}_{(j+1)}^{\alpha}
\widetilde{C}'_{\alpha-r-k}(\bar{q})\mathcal{N}^\alpha_1\f{\pa
f(\bar{q})}{\pa \bar{q}^l} - \mathcal{N}_{(l+1)}^{\alpha}
\widetilde{C}'_{\alpha-r-k}(\bar{q})\mathcal{N}^\alpha_1\f{\pa
f(\bar{q})}{\pa \bar{q}^j}  =0 \label{nec2} \ee Taking into
account  (\ref{implicit}), the above relation is fulfilled, for
any $j,l=1,...,m-1$. Indeed, taking the derivatives with respect
to $\bar{q}^j$ and also to $\bar{q}^l$, in the equation
(\ref{implicit}), we obtain, respectively, \ba \sum
^{b}_{\alpha=r+k+1}\mathcal{N}^\alpha_1\widetilde{C}'_{\alpha-r-k}(\bar{q})
\left[\mathcal{N}^\alpha_1\f{\pa f(\bar{q})}{\pa \bar{q}^j}
+\mathcal{N}_{(j+1)}^{\alpha}\right]&=&0\nonumber\\
\sum
^{b}_{\alpha=r+k+1}\mathcal{N}^\alpha_1\widetilde{C}'_{\alpha-r-k}(\bar{q})
\left[\mathcal{N}^\alpha_1\f{\pa f(\bar{q})}{\pa \bar{q}^l}
+\mathcal{N}_{(l+1)}^{\alpha}\right]&=&0\label{nec3} \ea

Now we multiply  in (\ref{nec3}) the first equation  with $\f{\pa
f(\bar{q})}{\pa \bar{q}^l}$, the second  equation with $-\f{\pa
f(\bar{q})}{\pa \bar{q}^j}$ and we add the resulting  equations to
obtain the equation (\ref{nec2}).

Thus, we proved  the existence of  a function
$\bar{E}_{0_\omega}(\bar{q},\dot{\bar{q}})$ such that
(\ref{dissip'}) is fulfilled, with the dissipative 1-form given by
(\ref{Dbar}). Therefore, the Birkhoffian (\ref{birred}) is
dissipative.

\noindent \noindent
 For the Birkhoffian (\ref{birred}),  the determinant in
 (\ref{regular}) becomes
\be
\textrm{det}\left[\f{\pa \bar{Q}_j}{\pa \ddot{\bar{q}}^i}(\bar{q},\, \dot{\bar{q}},\, 
\ddot{\bar{q}})\right]_{i,j=1,...,m-1}=\textrm{det}\left[
\sum^{r+k}_{a=r+1}\mathcal{N}^a_{(j+1)}\mathcal{N}^ {a}_{(i+1)}
\widetilde{L}_{a-r}\left(\dot{\bar{q}}\right) \right]\label{d}\ee
If the determinant in (\ref{d}) is different from zero, then the
Birkhoffian (\ref{birred}) is regular. $\quad\blacksquare$

 \textit{If
 there  exists in the network  $m_2<m$  loops which contain only resistors,
 all the other loops containing at least an inductor, we can
 \textbf{regularize} the Birkhoffian (\ref{bir}) via \textbf{reduction of the
 configuration space}. The reduced configuration space $\hat{M}_c$
 of dimension $m-m_2$, is a linear or a nonlinear subspace of $M_c$,
depending on whether the resistors  are linear or nonlinear. We
claim that  the Birkhoffian $\hat{\omega}_c$ on the reduced
configuration space $\hat{M}_c$ is still a \textbf{dissipative
Birkhoffian}. Under certain conditions on the functions $L_a$,
$a=1,...,k,$ which characterize the inductors, the reduced
Birkhoffian $\hat{\omega}_c$ will be a \textbf{regular
Birkhoffian}}.

Without loss of generality, we may assume that we have one loop in
the network that contains only resistors and in the coordinate
system that we have chosen, the constants read as \be
\mathcal{N}_1^{a}=0, \, a=r+1,..., r+k, \quad
\mathcal{N}_1^{\alpha}=0, \,\alpha=r+k+1,...,b\label{rlin2}\ee

I) Let us first consider the case that the resistors in this loop
are linear resistors, that is, described by (\ref{15rez}), all the
other electrical devices in the network being nonlinear. This
means that we have \be \mathcal{N}_1^{\Gamma}\neq 0, \,
\Gamma=1,...,r_{lin}, \quad \mathcal{N}_1^{\Gamma}=0, \,
\Gamma=r_{lin}+1,..., r\label{rlin1}\ee
 where $r_{lin}$ is the
number of linear resistors in the network. \\
In this case, the expressions of the Birkhoffian components
(\ref{bir}), with (\ref{bir2}), (\ref{bir3}), (\ref{bir1}), are
given by,  $j=2,...,m$, \ba Q_1(q,\dot{q},\ddot{q})&=&\sum
^{r_{lin}}_{\Gamma=1} \sum^{m}_{l=1}\mathcal{N}_1^{\Gamma} \textrm
{\scriptsize R}_{\Gamma}\mathcal{N}^ {\Gamma}_{l}\dot{q}^l
\nonumber\\
Q_j(q,\dot{q},\ddot{q})&=&\sum^{m}_{i=2}
\sum^{r+k}_{a=r+1}\mathcal{N}^a_{j}\mathcal{N}^
{a}_{i}\widetilde{L}_{a-r}\left(\dot{q}\right)\ddot{q}^i+\sum
^{r_{lin}}_{\Gamma=1} \sum^{m}_{l=1}\mathcal{N}_j^{\Gamma} \textrm
{\scriptsize R}_{\Gamma}\mathcal{N}^ {\Gamma}_{l}\dot{q}^l
+\nonumber\\
&&\sum ^{r}_{\Gamma=r_{lin}+1}\mathcal{N}_j^{\Gamma}
\widetilde{R}_{\Gamma}\left(\dot{q}\right)+\sum
^{b}_{\alpha=r+k+1}\mathcal{N}_j^{\alpha}\widetilde{C}_{\alpha-r-k}(q)
\label{resistorlinear} \ea We note that according to
(\ref{rlin2}), (\ref{rlin1}), $\dot{q}^1$ does  not appear in any
function $\widetilde{R}_{\Gamma}(\dot{q})$,
$\widetilde{L}_{a-r}(\dot{q})$, the terms
$\widetilde{L}_{a-r}(\dot{q})\ddot{q}^1$ do not appear in any of
the Birkhoffian components $Q_2(q,\dot{q}, \ddot{q}),..., Q_m(q,
\dot{q}, \ddot{q})$ as well $q^1$ does not
appear in any function $\widetilde{C}_{\alpha-r-k}(q)$.\\
 We define the $(m-1)$-dimensional
linear space $\hat{M}_{c}\subset M_c$  by \be \hat{M}_{c}=\{q\in
M_c\,\, | \, \, \sum ^{r_{lin}}_{\Gamma=1}
\sum^{m}_{l=1}\mathcal{N}_1^{\Gamma} \textrm {\scriptsize
R}_{\Gamma}\mathcal{N}^ {\Gamma}_{l}q^l+c_1=0\}\label{mbar} \ee
with $c_1$ a real constant.\\
We take $\hat{q}^1:=q^2$,..., $\hat{q}^{m-1}:=q^m$, as local
coordinates on the reduced configuration space $\hat{M}_c$. Then,
making use of  (\ref{mbar}) and  of the fact that
$\mathcal{N}_1^{\Gamma}\neq 0$ and $\textrm {\scriptsize
R}_{\Gamma}>0$, for any $\Gamma=1,...,r_{lin}$, we can express
$\dot{q}^1$  as a linear combination of
$\dot{\hat{q}}^1,...,\dot{\hat{q}}^{m-1} $, denoted
$g(\dot{\hat{q}})$, such that \be \sum ^{r_{lin}}_{\Gamma=1}
\mathcal{N}_1^{\Gamma} \textrm {\scriptsize
R}_{\Gamma}\left[\mathcal{N}^
{\Gamma}_{1}g(\dot{\hat{q}})+\sum^{m-1}_{l=1}\mathcal{N}^
{\Gamma}_{(l+1)}\dot{\hat{q}}^l\right]=0\label{g}\ee Thus, the
reduced Birkhoffian has the form
$\hat{\omega}_c=\sum^{m-1}_{j=1}\hat{Q}_jd\hat{q}^j$,
  \be
\hat{Q}_j(\hat{q}, \dot{\hat{q}},
\ddot{\hat{q}})=\hat{F}_j(\dot{\hat{q}})\ddot{\hat{q}}+
\hat{H}_j(\dot{\hat{q}})+\hat{G}_j(\hat{q}), \quad \textrm{where}
\label{birred2} \ee
\be 
\hat{F}_j(\dot{\hat{q}})\ddot{\hat{q}}:= \sum^{m-1}_{i=1}
\sum^{r+k}_{a=r+1}\mathcal{N}^a_{(j+1)}\mathcal{N}^ {a}_{(i+1)}
L_{a-r}\left(\sum^{m-1}_{l=1} \mathcal{N}^a_{(l+1)}
\dot{\hat{q}}^l\right) \ddot{\hat{q}}^i\ee \ba
\hat{H}_j(\dot{\hat{q}}):&=&\sum ^{r_{lin}}_{\Gamma=1}
\mathcal{N}_{(j+1)}^{\Gamma} \textrm {\scriptsize
R}_{\Gamma}\left[\mathcal{N}^
{\Gamma}_{1}g(\dot{\hat{q}})+\sum^{m-1}_{l=1}\mathcal{N}^
{\Gamma}_{(l+1)}\dot{\hat{q}}^l\right]+ \nonumber\\
&&\sum ^{r}_{\Gamma=r_{lin}+1}\mathcal{N}_{(j+1)}^{\Gamma}
R_{\Gamma}\left(\sum^{m-1}_{l=1} \mathcal{N}^\Gamma_{(l+1)}
\dot{\hat{q}}^l\right) \label{H} \ea \be  \hat{G}_j(\hat{q}):=\sum
^{b}_{\alpha=r+k+1} \mathcal{N}_{(j+1)}^{\alpha} C_{\alpha-r-k}
\left(\sum^{m-1}_{l=1} \mathcal{N}^\alpha_{(l+1)}
\hat{q}^l+const\right)\label{bir66'} \ee  The Birkhoffian given by
(\ref{birred2}) is still dissipative. We will see that there
exists  a function $\hat{E}_{0_\omega}(\hat{q},\dot{\hat{q}})$
such that \be \sum^{m-1}_{j=1}\hat{Q}_j(\hat{q}\, \,
\dot{\hat{q}},\,
\ddot{\hat{q}})\dot{\hat{q}}^j=\sum^{m-1}_{j=1}\left[\f{\pa
\hat{E}_{0_\omega}}{\pa \hat{q}^j}\dot{\hat{q}}^j+\f{\pa
\hat{E}_{0_\omega}}{\pa \dot{\hat{q}}^j}\ddot{\hat{q}}^j+
\hat{D}_{j}(\hat{q},\dot{\hat{q}})\dot{\hat{q}}^j\right]
\label{dissip1}\ee where
$\hat{D}=\sum_{j=1}^{m-1}\hat{D}_{j}(\hat{q},\dot{\hat{q}})d{\hat{q}}^j$
is a dissipative 1-form on $T\hat{M_c}$.\\
We consider the following   1-form on $T\hat{M}_c$ \be
\hat{D}=\sum_{j=1}^{m-1}\hat{H}_j(\dot{\hat{q}})
d\hat{q}^j\label{disipbar} \ee Let us check that the vertical
1-form (\ref{disipbar}) is dissipative, that is,
  \be
\sum_{j=1}^{m-1}\hat{H}_j(\dot{\hat{q}})\dot{\hat{q}}^j>
0\label{dhat} \ee From (\ref{H}), the left  side of (\ref{dhat})
writes as the sum $\mathcal{S}_1+\mathcal{S}_2$, where \be
\mathcal{S}_1=\sum_{j=1}^{m-1}\sum ^{r_{lin}}_{\Gamma=1}
\mathcal{N}_{(j+1)}^{\Gamma} \textrm {\scriptsize
R}_{\Gamma}\left[\mathcal{N}^
{\Gamma}_{1}g(\dot{\hat{q}})+\sum^{m-1}_{l=1}\mathcal{N}^
{\Gamma}_{(l+1)}\dot{\hat{q}}^l\right]
\dot{\hat{q}}^j\label{s1}\ee \be \mathcal{S}_2=\sum
^{r}_{\Gamma=r_{lin}+1} \left[
R_{\Gamma}\left(\sum_{j=1}^{m-1}\mathcal{N}_{(j+1)}^{\Gamma}\dot{\hat{q}}^j\right)\right]\left(
\sum_{j=1}^{m-1}\mathcal{N}_{(j+1)}^{\Gamma}\dot{\hat{q}}^j\right)\label{s2}\ee
We  now multiply  the equation  (\ref{g}) by the function
$g(\dot{\hat{q}})$. Using the resulting equation  we can write the
sum in (\ref{s1}) in the form \be \mathcal{S}_1=\sum
^{r_{lin}}_{\Gamma=1}  \textrm {\scriptsize
R}_{\Gamma}\left[\mathcal{N}^
{\Gamma}_{1}g(\dot{\hat{q}})+\sum^{m-1}_{l=1}\mathcal{N}^
{\Gamma}_{(l+1)}\dot{\hat{q}}^l\right]^2\ee Since $\textrm
{\scriptsize R}_{\Gamma}>0$, $\Gamma=1,...,r_{lin}$, the sum
$\mathcal{S}_1$  is strictly positive.\\
Because all  nonlinear resistors considered satisfy the condition
(\ref{assumption}),  the  sum $\mathcal{S}_2$ in (\ref{s2}) is
strictly positive as well. Therefore, the inequality (\ref{dhat})
is fulfilled and the vertical 1-form in (\ref{disipbar}) is
dissipative.\\
We now look for a function
$\hat{E}_{0_\omega}(\hat{q},\dot{\hat{q}})$ such that \be
\sum^{m-1}_{j=1}\left[
\hat{F}_j(\dot{\hat{q}})\ddot{\hat{q}}+\hat{G}_j(\hat{q})\right]
\dot{\hat{q}}^j=\sum^{m-1}_{j=1}\left[\f{\pa
\hat{E}_{0_\omega}}{\pa \hat{q}^j}\dot{\hat{q}}^j+\f{\pa
\hat{E}_{0_\omega}}{\pa
\dot{\hat{q}}^j}\ddot{\hat{q}}^j\right]\label{dissip'''} \ee
Because of the special form of the terms on the left  side of (\ref{dissip'''}), we can look for the required function 
$\hat{E}_{0_\omega}(\hat{q},\dot{\hat{q}})$ as a sum of a function
only depending on $\hat{q}$, and a function only depending on
$\dot{\hat{q}}$.
From the theory of 
total differentials, a necessary condition for the existence of
such functions is the fulfilment of the following relations \be
\left\{
\begin{array}{ll}
\f {\partial \hat{\mathcal{F}}_j(\dot{\hat{q}})}{\partial
\dot{\hat{q}}^l }-\f {\partial \hat{\mathcal{F}}_l
(\dot{\hat{q}})}{\partial
\dot{\hat{q}}^j}=0\\
\\
\f {\partial \hat{G}_j(\hat{q})}{\partial \hat{q}^l}-\f {\partial
\hat{G}_l(\hat{q})}{\partial \hat{q}^j}=0
\end{array}
\right.\label{bir55} \ee for any $j,l=1,...,m-1$, where \be
\hat{\mathcal{F}}_j(\dot{\hat{q}})=\sum^{m-1}_{i=1}\sum^{r+k}_{a=r+1}\mathcal{N}^{\,
a}_{(j+1)}\mathcal{N}^ {a}_{(i+1)} L_{a-r} \left(\sum^{m-1}_{l=1}
\mathcal{N}^a_{(l+1)} \dot{\hat{q}}^l\right) \dot{\hat{q}}^i
\label{bir66} \ee  From (\ref{bir66}), (\ref{bir66'}) we get \ba
\hspace{-0.7cm}\f {\partial
\hat{\mathcal{F}}_j(\dot{\hat{q}})}{\partial \dot{\hat{q}}^l}=
\sum^{r+k}_{a=r+1}\mathcal{N}^a_{(j+1)}\mathcal{N}^{a}_{(l+1)}\widetilde{L}_{a-r}(\dot{\hat{q}})+
\sum^{m-1}_{i=1}\sum^{r+k}_{a=r+1}\mathcal{N}^a_{(j+1)}\mathcal{N}^
{a}_{(i+1)}\mathcal{N}^
{a}_{(l+1)}\widetilde{L}_{a-r}'(\dot{\hat{q}})\dot{\hat{q}}^i\nonumber\\\label{bir77}
\ea \be \hspace{-0.7cm}\f {\partial \hat{G}_j(\hat{q})}{\partial
\hat{q}^l}= \sum
^{b}_{\alpha=r+k+1}\mathcal{N}^{\alpha}_{(j+1)}\mathcal{N}^{\alpha}_{(l+1)}\widetilde{C}_{\alpha-r-k}'(\hat{q})
\ee where  $\widetilde{L}'_{a-r}:=\f
{d\widetilde{L}_{a-r}(\eta)}{d\eta}$,
$\widetilde{C}'_{\alpha-r-k}:=\f
{d\widetilde{C}_{\alpha-r-k}(\eta)}{d\eta}$.
We can easily check that the  equations (\ref{bir55}) are
fulfilled.\\
Thus, we proved  the existence of  a function
$\hat{E}_{0_\omega}(\hat{q},\dot{\hat{q}})$ such that
(\ref{dissip1}) is fulfilled, with the dissipative 1-form given by
(\ref{disipbar}), that is, the Birkhoffian (\ref{birred2}) is
dissipative.

 For the Birkhoffian (\ref{birred2}),  the determinant in
 (\ref{regular}) becomes
\be
\textrm{det}\left[\f{\pa \hat{Q}_j}{\pa \ddot{\hat{q}}^i}(\hat{q},\, \dot{\hat{q}},\, 
\ddot{\hat{q}})\right]_{i,j=1,...,m-1}=\textrm{det}\left[
\sum^{r+k}_{a=r+1}\mathcal{N}^a_{(j+1)}\mathcal{N}^ {a}_{(i+1)}
\widetilde{L}_{a-r}\left(\dot{\bar{q}}\right) \right]\label{d'}\ee
If the determinant in (\ref{d'}) is different from zero, then the
Birkhoffian (\ref{birred2}) is regular. $\quad\blacksquare$

 II) Let us now consider the case that the resistors in the
loop formed only by resistors are  nonlinear devices too, that is,
they are described by (\ref{2rezistor}), with the assumption
(\ref{assumption}). Now, instead of (\ref{rlin1}) we have $
\mathcal{N}_1^{\Gamma}\neq 0, \, \Gamma=1,...,r$.\\
 The component
$Q_1(q,\dot{q},\ddot{q})$ of the Birkhoffian takes the form \be
Q_1(q,\dot{q},\ddot{q})=\sum ^{r}_{\Gamma=1}
\mathcal{N}_1^{\Gamma} \widetilde{R}_{\Gamma}(\dot{q})\label{dotq}
\ee and
 the other components $Q_2(q,\dot{q},\ddot{q})$,...,
$Q_m(q,\dot{q},\ddot{q})$ are the same as in
(\ref{resistorlinear}) with the terms following  $r_{lin}=0$
absent. According to (\ref{rlin2}), $\dot{q}^1$ does not appear in
any function $\widetilde{L}_{a-r}(\dot{q})$, the terms
$\widetilde{L}_{a-r}(\dot{q})\ddot{q}^1$ do not appear in any of
the Birkhoffian components $Q_2(q,\dot{q}, \ddot{q}),..., Q_m(q,
\dot{q}, \ddot{q})$ as well $q^1$ does not appear in any function
$\widetilde{C}_{\alpha-r-k}(q)$.\\
 Using (\ref{dotq}), we intend to  define the
$(m-1)$-dimensional configuration space  $\hat{M}_{c}\subset M_c$.
The relation (\ref{dotq}) is a nonlinear velocity constraint,
which in general is  a nonholonomic constraint. Nevertheless,
because of this constraint imposed on the system, the equations
which describe the dynamics are\\
$Q_2(q,\dot{q},\ddot{q})=0$,..., $Q_m(q,\dot{q},\ddot{q})=0$.\\
Taking $\hat{q}^1:=q^2$,..., $\hat{q}^{m-1}:=q^m$, a coordinate
system on the reduced configuration space $\hat{M}_c$, the
Birkhoffian has in this case the form
$\hat{\omega}_c=\sum^{m-1}_{j=1}\hat{Q}_jd\hat{q}^j$,  \ba
\hat{Q}_j(\hat{q}, \dot{\hat{q}}, \ddot{\hat{q}})&=&
\sum^{m-1}_{i=1}
\sum^{r+k}_{a=r+1}\mathcal{N}^a_{(j+1)}\mathcal{N}^ {a}_{(i+1)}
L_{a-r}\left(\sum^{m-1}_{l=1} \mathcal{N}^a_{(j+1)}
\dot{\hat{q}}^l\right)
\ddot{\hat{q}}^i+\nonumber\\
&&\sum ^{r}_{\Gamma=1}\mathcal{N}_{(j+1)}^{\Gamma}
R_{\Gamma}\left(\mathcal{N}^
{\Gamma}_{1}h(\dot{\hat{q}}^1,...,\dot{\hat{q}}^m)+\sum^{m-1}_{l=1}\mathcal{N}^
{\Gamma}_{(l+1)}\dot{\hat{q}}^l\right)\nonumber\\
&& +\sum ^{b}_{\alpha=r+k+1} \mathcal{N}_{(j+1)}^{\alpha}
C_{\alpha-r-k} \left(\sum^{m-1}_{l=1} \mathcal{N}^\alpha_{(l+1)}
\hat{q}^l\right)
 \label{birrednonlinear}
\ea where $h:U\subset\mathbf{R}^{m-1}\longrightarrow \mathbf{R}$
is the unique function such that $h(\dot{\hat{q}}_0)=\dot{q}^1_0$,
$\dot{q}^1_0\in \mathbf{R}$, and \be \sum ^{r}_{\Gamma=1}
\mathcal{N}_1^{\Gamma} R_{\Gamma}
\left(\mathcal{N}^\Gamma_1h(\dot{\hat{q}}^1,...,\dot{\hat{q}}^{m-1})+\sum^{m-1}_{l=1}
\mathcal{N}^\Gamma_{(l+1)}
\dot{\hat{q}}^l\right)=0\label{implicit1} \ee
 for all $\dot{\hat{q}}=(\dot{\hat{q}}^1,..., \dot{\hat{q}}^{m-1})\in U$, with
$U$ a neighborhood of
$\dot{\hat{q}}_0=(\dot{\hat{q}}_0^1,...,\dot{\hat{q}}_0^{m-1})$.\\
One can prove using the same ideas as in the previous case, that
the Birkhoffian given by (\ref{birrednonlinear}) is still
dissipative. \\
If the determinant  $\textrm{det} \left[
\sum^{r+k}_{a=r+1}\mathcal{N}^a_{(j+1)}\mathcal{N}^ {a}_{(i+1)}
\widetilde{L}_{a-r}\left(\dot{\hat{q}}\right)
\right]_{i,j=1,...,m-1}\neq 0$, then, the Birkhoffian
(\ref{birrednonlinear}) is regular.
 $\quad\blacksquare$

 \textit{If
 there  exists in the network  $m_3<m$  loops which contain only resistors and capacitors,
 all the other loops containing at least an inductor, we can
 \textbf{regularize} the Birkhoffian (\ref{bir})
by introducing  into each of these $m_3$ loops,  an inductor in
series, with the  inductance functions
$\mathcal{L}_{a'}:\mathbf{R}\longrightarrow \mathbf{R}\backslash
\{0\}$, $a'=1,...,m_3$, having very small values. The
configuration space remains $M_{c}$
 of dimension $m$. under certain conditions on the functions
 $\mathcal{L}_{a'}$ and on the functions $L_{a}$,
$a=1,...,k$, which characterize the others inductors, the
Birkhoffian $\omega^{ext}_{c}$ on  $M_c$ will
 be a \textbf{dissipative regular
Birkhoffian}}.

 Without loss of generality, we may assume that we have one
loop in the network  containing  only resistors and capacitors and
in the chosen coordinate system
 \be
\mathcal{N}_1^{a}=0, \, a=r+1,..., r+k \ee \vspace{-0.5cm}
 The expressions
(\ref{bir}) of the Birkhoffian components become, $j=2,...,m.$,
 \ba
&&\hspace{-1cm}Q_1(q,\dot{q},\ddot{q})= \sum ^{r}_{\Gamma=1}
\mathcal{N}_1^{\Gamma} \widetilde{R}_{\Gamma}(\dot{q})+\sum
^{b}_{\alpha=r+k+1}\mathcal{N}_1^{\alpha}\widetilde{C}_{\alpha-r-k}(q)
\nonumber\\
&&\hspace{-1cm}Q_j(q,\dot{q},\ddot{q})=\sum^{m}_{i=2}
\sum^{r+k}_{a=r+1}\mathcal{N}^a_{j}\mathcal{N}^
{a}_{i}\widetilde{L}_{a-r}\left(\dot{q}\right)\ddot{q}^i+ \sum
^{r}_{\Gamma=1}\mathcal{N}_j^{\Gamma}
\widetilde{R}_{\Gamma}\left(\dot{q}\right) +\sum
^{b}_{\alpha=r+k+1}\mathcal{N}_j^{\alpha}\widetilde{C}_{\alpha-r-k}(q)
\nonumber\\
&&\label{rescap1} \ea
 We introduce into this loop an inductor in
series, described by the following relation between the current
and the  voltage \be
v=\mathcal{L}_{1}({\textsc{i}})\f{d{\textsc{i}}}{dt} \ee
$\mathcal{L}_{1}:\mathbf{R}\longrightarrow \mathbf{R}\backslash
\{0\}$ being smooth invertible function. After introducing the
inductor, the number of branches of the graph associated to the
circuit increases by one, that is, there will be $b+1$ branches,
and the number of nodes increases by one as well, that is, $n$
becomes $n+1$. Still the cardinality of a selection of loops which
cover the whole graph remains $m$. The configuration space is the
same $M_c$, with dimension $m$. The corresponding Birkhoffian,
denoted by $\omega^{ext}_c$, has  the component
$Q_1(q,\dot{q},\ddot{q})$ given by,
 \be
Q_1(q,\dot{q},\ddot{q})=
\widetilde{\mathcal{L}}(\dot{q}^1)\ddot{q}^1+\sum ^{r}_{\Gamma=1}
\mathcal{N}_1^{\Gamma} \widetilde{R}_{\Gamma}(\dot{q})+\sum
^{b}_{\alpha=r+k+1}\mathcal{N}_1^{\alpha}\widetilde{C}_{\alpha-r-k}(q)\ee
the others $Q_2(q,\dot{q},\ddot{q})$,...,
$Q_m(q,\dot{q},\ddot{q})$ have the same form as in
(\ref{rescap1}). This Birkhoffian is dissipative, the expression
of the function $E_{0_\omega}$ on $TM_c$ is (\ref{energie}) plus
the term $\int
\widetilde{\mathcal{L}}(\dot{q}^1)\dot{q}^1d\dot{q}^1$. The
dissipative 1-form has  the form (\ref{linearvert}).\\
If $\widetilde{\mathcal{L}}(\dot{q}^1)\textrm{det} \left[
\sum^{r+k}_{a=r+1}\mathcal{N}^a_{(j+1)}\mathcal{N}^ {a}_{(i+1)}
\widetilde{L}_{a-r}\left(\dot{\hat{q}}\right)
\right]_{i,j=1,...,m-1}\neq 0$,  this Birkhoffian
 is regular. $\quad\blacksquare$
\subsection{Voltage controlled resistors}
Let us now consider the nonlinear resistors for which the
constitutive relations are given by \be
\textsc{i}_{\Gamma}=\mathfrak{R}_{\Gamma}(v_{\Gamma}), \quad
\Gamma= 1,...,r \label{2vcr} \ee where
$\mathfrak{R}_{\Gamma}:\mathbf{R}\longrightarrow \mathbf{R}$ are
smooth functions. In order to obtain  a dissipative Birkhoffian,
we also assume that, for all $x \neq 0$, \be
\mathfrak{R}_{\Gamma}(x)x>0, \quad \forall
\Gamma=1,...,r\label{assumptionvcr} \ee that is, for each
nonlinear resistor, the graph of the function
$\mathfrak{R}_{\Gamma}$ lies in the union of the  first and the third quadrant.\\
Taking into account (\ref{2}), (\ref{4'}), (\ref{4}),
(\ref{2vcr}),  the equations (\ref{5}), (\ref{6})   governing the
circuit have the form \be \left\{
\begin{array}{ll}
B^T
\left(
\begin{array}{c}
\textsc{i}_{\Gamma}\\
\textsc{i}_{a} \\
\dot{\textsc{q}}_{\alpha}
\end{array}
\right)=0\\
\\
\left\{
\begin{array}{c}
A^T
\left(
\begin{array}{c}
v_{\Gamma}\\
L_{a}(\textsc{i}_a)\,
\dot{\textsc{i}_a} \\
C_\alpha(\textsc{q}_{\alpha})
\end{array}
\right)=0\\
\\
\textsc{i}_{\Gamma}=\mathfrak{R}_{\Gamma}(v_\Gamma)
\end{array}
\right.
\end{array}
\right. \label{7vcr} \ee
 \textit{Using the first set of equations
(\ref{7vcr}),
 we define a  family of
  $m$-dimensional affine-linear configuration spaces
  $M_c\subset\mathbf{R}^b$ parameterized by a constant vector
  $c$ in $\mathbf{R}^n$. This vector is related to the initial
  values of the $\textsc{q}$-variables at some instant of time.
A Birkhoffian $\omega_c$ on the configuration space $M_c$ arises
from the second set of equations (\ref{7vcr}). Thus,
($M_c,\omega_c$) will be a family of Birkhoff systems that
describe the considered RLC network.}

The first  set of equations (\ref{7vcr}) is the same as the first
set of equations  (\ref{7}). Thus, using them we can define the
family $M_c$ of $m$-dimensional affine-linear configuration spaces
(\ref{8}). For a coordinate system $q=(q^1,...,q^m)$ on $M_c$, the
relations between the $x$-coordinates (\ref{cs}) and the
$q$-coordinates  are given by (\ref{10'}). The matrix of constants
$\mathcal{N}$ satisfies (\ref{N}).
 Taking into account (\ref{ker}),  we define the
Birkhoffian $\omega_c$ of $M_c $ such that the differential system
(\ref{difsistem}) is the linear combination of the second set of
equations in (\ref{7vcr}) obtained by replacing $A^T$ with the
matrix $\mathcal{N}^T$. In terms of the $q$-coordinates chosen
before, the components  $Q_j(q,\dot{q},\ddot{q})$ of the
Birkhoffian have the implicit form
 \be \left\{
\begin{array}{ll}
Q_j(q, \, \dot{q},\,\ddot{q})=
F_j(\dot{q})\ddot{q}+G_j(q)+
\sum ^{r}_{\Gamma=1}\mathcal{N}_j^{\Gamma}v_{\Gamma}, \quad j=1,...,m\\
\\
\sum ^{m}_{j=1}\mathcal{N}_j^{\Gamma}\dot{q}^j=\mathfrak{R}_{\Gamma}
(v_\Gamma)
\end{array}
\right. \label{birvcr} \ee where the functions
$F_j(\dot{q})\ddot{q}$, $G_j(q)$ are given by (\ref{bir2}),
(\ref{bir1}).

\textit{  The Birkhoffian (\ref{birvcr}) is \textbf{dissipative}
}.

Indeed, with the  function $E_{0_\omega}(q,\dot{q})$ given by
(\ref{energie}), the identity (\ref{dissip}) becomes for the
Birkhoffian (\ref{birvcr}), \be \sum_{j=1}^{m}Q_j(q, \,
\dot{q},\,\ddot{q})\dot{q}^j= \sum_{j=1}^m\left[\f{\pa
E_{0_\omega}}{\pa q^j}\dot{q}^j+\f{\pa E_{0_\omega}}{\pa
\dot{q}^j}\ddot{q}^j+ \sum
^{r}_{\Gamma=1}\mathcal{N}_j^{\Gamma}v_{\Gamma}\dot{q}^j\right]
\ee It remains to show that the vertical 1-form defined implicitly
 by \be \left\{
\begin{array}{ll}
D_j(q, \, \dot{q})=
\sum ^{r}_{\Gamma=1}\mathcal{N}_j^{\Gamma}v_{\Gamma}, \quad j=1,...,m\\
\\
\sum
^{m}_{j=1}\mathcal{N}_j^{\Gamma}\dot{q}^j=\mathfrak{R}_{\Gamma}
(v_\Gamma)
\end{array}
\right. \label{birvcr'} \ee is dissipative. From the second set of
relations in (\ref{birvcr'}), we have \be \sum_{j=1}^m\sum
^{r}_{\Gamma=1}\mathcal{N}_j^{\Gamma}v_{\Gamma}\dot{q}^j= \sum
^{r}_{\Gamma=1}\mathfrak{R}_{\Gamma} (v_\Gamma)v_{\Gamma} \ee
Therefore, the inequality (\ref{store}) reads as $\sum
^{r}_{\Gamma=1}\mathfrak{R}_{\Gamma} (v_\Gamma)v_{\Gamma} >0.$ The
last inequality is satisfied  in view of the assumption
(\ref{assumptionvcr}).  $\quad\blacksquare$

\textit{As in the case of current controlled sources, the
Birkhoffian (\ref{birvcr}) is \textbf{\textit{never} regular} if
the network contains closed loops formed only by capacitors, or
resistors, or both of them.}

\section{Examples}

Example 1):   This example is based on the following oriented
connected graph

\begin{center}
\scalebox{0.40}{\includegraphics{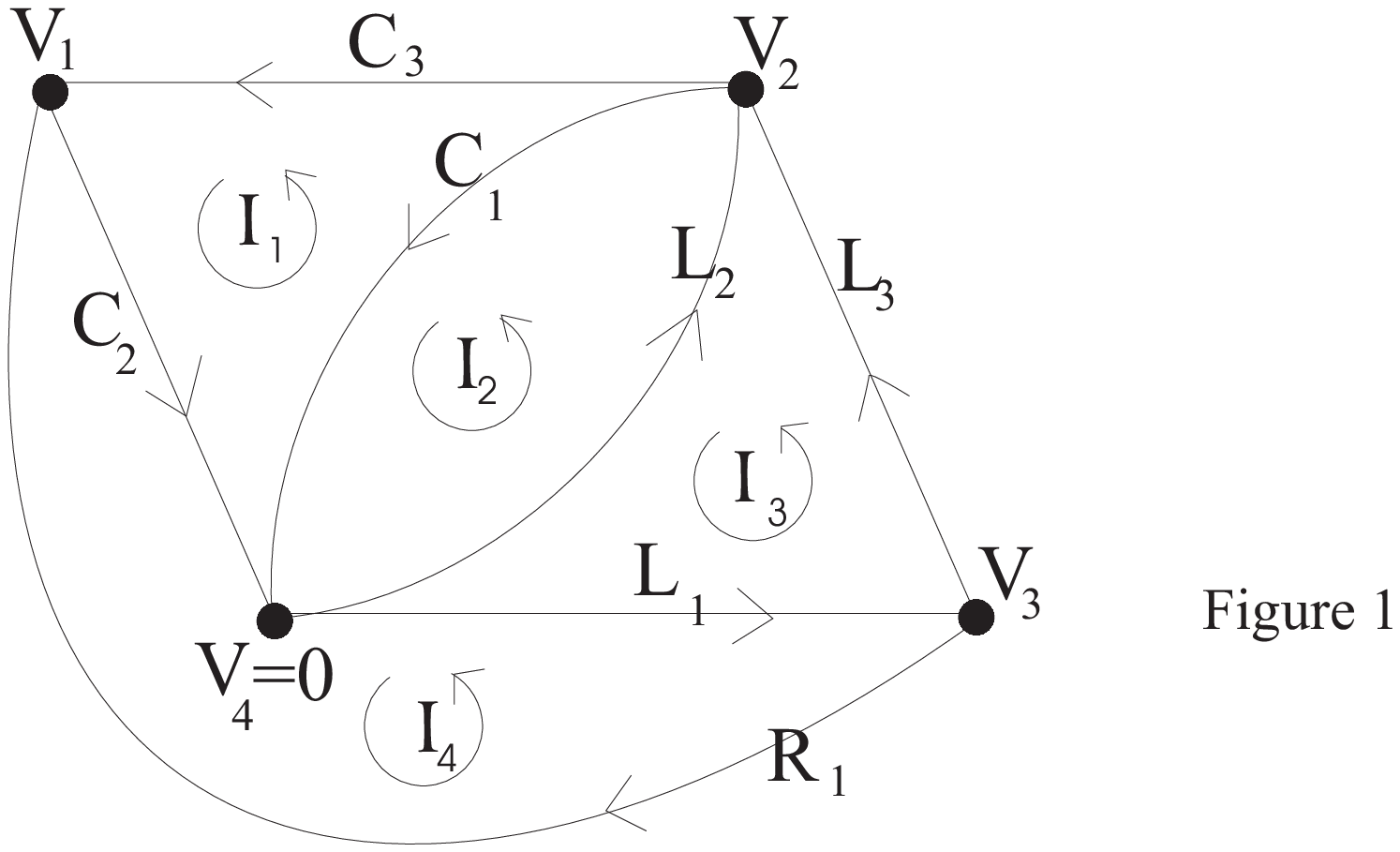}}
\end{center}

\noindent We have $r=1$, $k=3$, $p=3$, $n=3$, $m=4$ , $b=7$. We
choose the reference node to be $V_4$ and the current directions
as indicated in Figure 1. We cover the associated graph with the
loops $I_1, \, I_2,\, I_3,\, I_4$. Let $V=(V_1,V_2,V_3)\in
\mathbf{R}^3$ be the vector of node voltage values,
$\textsc{i}=(\textsc{i}_{[\Gamma]},
\textsc{i}_{(a)},\textsc{i}_\alpha)\in \mathbf{R}^1\times
\mathbf{R}^3\times \mathbf{R}^3$ be
 the vector of branch current  values  and
 $v=(v_{[\Gamma]},v_({a}),v_{\alpha})\in \mathbf{R}^1\times \mathbf{R}^3\times\mathbf{R}^3$ be the vector of branch  voltage
 values.\\
The branches in Figure 1 are  labelled as follows: the first
branch is the resistive branch $\textsc{r}_1$, the second, the
third and the fourth branch are the inductive branches
$\textsc{l}_1$, $\textsc{l}_2$, $\textsc{l}_3$ and the last three
branches are the capacitor branches $\textsc{c}_1$,
$\textsc{c}_2$, $\textsc{c}_3$. The incidence and loop matrices,
$B\in \mathfrak{M}_{73}(\mathbf{R})$ and
 $A\in \mathfrak{M}_{74}(\mathbf{R})$, write as \be
B=\left(
\begin{array}{ccc}
-1& 0& 1\\
0& 0&-1\\
 0& -1& 0 \\
0&-1& 1\\
 0&1& 0
 \\
 1& 0& 0 \\
-1& 1& 0
\end{array}
\right),\quad \quad
 A=\left(
\begin{array}{cccc}
 0& 0 & 0 & -1\\
 0& 0 & 1 & -1 \\
0& 1 & -1 & 0 \\
 0& 0 & 1 & 0\\
 -1& 1 &0 & 0 \\
 1& 0 & 0 & -1\\
1& 0&  0& 0
\end{array}
\right) \label{exrlc31} \ee
One has  rank$(B)=3$,  rank$(A)=4$.\\
All the electrical devices are considered to be nonlinear and
described by the relations,
 (\ref{2}),  (\ref{4}), with  $C_1,C_2,C_3:\mathbf{R}\rightarrow \mathbf{R}\backslash \{0\}$,
  $L_1,L_2, L_3:\mathbf{R}\rightarrow \mathbf{R}\backslash \{0\}$ smooth invertible functions
  and by the relation (\ref{2rezistor}), with $R_{1}:\mathbf{R}\rightarrow \mathbf{R}$ smooth function,
  such that (\ref{assumption}) is satisfied, that is,
  for any $x\neq 0$
\be R_{1}(x)x>0\label{R1} \ee  The equations (\ref{7}) which
govern the network have the form \be \left\{
\begin{array}{llllllll}
-\textsc{i}_{[1]}+\dot{\textsc{q}}_{2}-\dot{\textsc{q}}_3=0\\
-\textsc{i}_{(2)}-\textsc{i}_{(3)}+\dot{\textsc{q}}_{1}+\dot{\textsc{q}}_3=0\\
\textsc{i}_{[1]}-\textsc{i}_{(1)}+\textsc{i}_{(3)}=0\\
\\
-C_1(\textsc{q}_1)+C_2(\textsc{q}_2)+C_3(\textsc{q}_3)=0\\
L_2(\textsc{i}_{(2)})\dot{\textsc{i}}_{(2)}+C_1(\textsc{q}_1)=0\\
L_1(\textsc{i}_{(1)})\dot{\textsc{i}}_{(1)}-
L_2(\textsc{i}_{(2)})\dot{\textsc{i}}_{(2)}+
L_3(\textsc{i}_{(3)})\dot{\textsc{i}}_{(3)}=0\\
-R_1(\textsc{i}_{[1]})-L_1(\textsc{i}_{(1)})\dot{\textsc{i}}_{(1)}-C_2(\textsc{q}_2)=0
\end{array}
\right.\label{exrlc32} \ee
The relations (\ref{var}), (\ref{cs})
read as follows for this example \be
\textsc{i}_{[1]}:=\dot{\textsc{q}}_{[1]},\quad
\textsc{i}_{(a)}:=\dot{\textsc{q}}_{(a)}, \quad
a=1,2,3\label{varrlc3} \ee \be x^1:=\textsc{q}_{[1]},\,
x^2:=\textsc{q}_{(1)},\,
x^3:=\textsc{q}_{(2)},\,x^{4}:=\textsc{q}_{(3)},\,
x^{5}:=\textsc{q}_{1},\,x^{6}:=\textsc{q}_{2},
x^{7}:=\textsc{q}_{3}\label{csrlc3} \ee
  Using the
first 3 equations of the system (\ref{exrlc32}) we define the
$4$-dimensional affine-linear configuration space $M_c$. In view
of the notations (\ref{varrlc3}), (\ref{csrlc3}), we integrate
these 3 equations and solving them  in terms of $4$ variables, we
obtain, for example, $ x^2=x^1+x^4+const,\,
x^5=x^{3}+x^{4}-x^{7}+const,\, x^6=x^1+x^{7}+const .$ Thus, a
coordinate system on $M_c$ is given by \be q^1:=x^7, q^2:=x^4,
q^3:=x^1,q^4:=x^3\label{qcoord2} \ee The matrix of constants
$\mathcal{N}=\left(\begin{array}{c}
\mathcal{N}^{\Gamma}_j\\
\mathcal{N}^{a}_{ j}\\
\mathcal{N}^{\alpha}_{j}
\end{array}\right)_{{\Gamma=1,2,a=3,4, \alpha=5,6,7 \atop
j=1,2,3,4}}$  from (\ref{10'}) is given by $ \mathcal{N}=\left(
\begin{array}{cccc}
 0& 0& 1& 0\\
 0& 1& 1& 0\\
 0& 0& 0& 1\\
 0& 1& 0& 0\\
 -1&1& 0& 1\\
 1& 0& 1& 0\\
 1& 0& 0& 0
\end{array}
\right)$.  Therefore, in terms of the $q$-coordinates
(\ref{qcoord2}),  we may define the Birkhoffian
$\omega_c=Q_1(q,\dot{q},\ddot{q})dq^1+Q_2(q,\dot{q},\ddot{q})dq^2+Q_3(q,\dot{q},\ddot{q})dq^3+Q_4(q,\dot{q},\ddot{q})dq^4$
of $M_c$ as in (\ref{bir}), with (\ref{bir2}), (\ref{bir3}),
(\ref{bir1}), that is, \ba
&&\hspace{-1cm}Q_1(q,\dot{q},\ddot{q})=-C_1(-q^1+q^2+q^4+const)+C_2(q^1+q^3+const)+
C_3(q^1)\nonumber\\
&&\hspace{-1cm}Q_2(q,\dot{q},\ddot{q})=\left[L_1(\dot{q}^2+\dot{q}^3)+
L_3(\dot{q}^2)\right]\ddot{q}^2+L_1(\dot{q}^2+\dot{q}^3)\ddot{q}^3+C_1(-q^1+q^2+q^4+const)
\nonumber\\
&&\hspace{-1cm}Q_3(q,\dot{q},\ddot{q})=L_1(\dot{q}^2+\dot{q}^3)\ddot{q}^2+
L_1(\dot{q}^2+\dot{q}^3)\ddot{q}^3+R_1(\dot{q}^3)+C_2(q^1+q^3+const)\nonumber\\
&&\hspace{-1cm}Q_4(q,\dot{q},\ddot{q})=L_2(\dot{q}^4)\ddot{q}^4+C_1(-q^1+q^2+q^4+const)
\label{exrlc36} \ea
 The Birkhoffian (\ref{exrlc36}) is \textbf{dissipative}
and \textbf{\textit{not} regular}.\\
 Indeed, there exists a smooth function
$E_{0_\omega}:TM\longrightarrow \mathbf{R}$ of the form
(\ref{energie}), that is,   \ba \hspace{-0.7cm}
E_{0_\omega}(q,\dot{q})&=& \int
\widetilde{L}_1(\dot{q}^2,\dot{q}^3)(\dot{q}^2+\dot{q}^3)(d\dot{q}^2+d\dot{q}^3)+\int
L_2(\dot{q}^4)\dot{q}^4d\dot{q}^4 +\int
L_3(\dot{q}^2)\dot{q}^2d\dot{q}^2
-\nonumber\\
&&
\int \int \widetilde{L}_1'(\dot{q}^2,\dot{q}^3)(\dot{q}^2+\dot{q}^3)d\dot{q}^2d\dot{q}^3-
\int \int \widetilde{L}_1(\dot{q}^2,\dot{q}^3)d\dot{q}^2d\dot{q}^3+\nonumber\\
&&
\int \widetilde{C}_1(q^1,q^2,q^4)(dq^1-dq^2+dq^4)+\int \widetilde{C}_2(q^1,q^3)(dq^1+dq^3)+
\nonumber\\
&&\int C_3(q^1)dq^1-\int \int
\widetilde{C}_1'(q^1,q^2,q^4)(-dq^1dq^2+dq^1dq^4-dq^2dq^4)-\nonumber\\
&&\int\int \widetilde{C}_2'(q^1,q^3)dq^1dq^3-\int\int\int
\widetilde{C}_1''(q^1,q^2,q^4)dq^1dq^2dq^4 \label{exrlc37} \ea
such that (\ref{dissip}) is satisfied with \be
D=R_1(\dot{q}^3)dq^3\label{d1} \ee Because the function $R_1$
satisfies (\ref{R1}), we obtain \be
D_j(q,\dot{q})\dot{q}^j=R_1(\dot{q}^3)\dot{q}^3>0 \ee that is,
(\ref{d1}) is indeed a dissipative vertical 1-form.\\
We are in the case where   the  circuit considered has one loop
which contains only capacitors. We note that for the Birkhoffian
(\ref{exrlc36}), the first row of the matrix $\left[\f{\pa
Q_j}{\pa \ddot{q}^i}\right]_{i,j=1,2,3,4}$ contains only zeros.
Therefore, $\textrm{det}\left[\f{\pa Q_j}{\pa
\ddot{q}^i}\right]_{i,j=1,2,3,4}=0$ and the Birkhoffian
(\ref{exrlc36}) is not regular.\\
Using the first relation in (\ref{exrlc36}), we now define a
3-dimensional $\bar{M}_{c}\subset M_c$ by \ba
\bar{M}_c=\{q=(q^1,q^2,q^3,q^4)\in
M_c/&&\,C_1(q^1-q^2+q^4+const)+\nonumber\\
&&C_2(q^1+q^3+const)+ C_3(q^1)=0 \}\label{barMrlc3} \ea By the
implicit function theorem, we obtain a  local coordinate system on
the reduced configuration space $\bar{M}_c$. Taking
$\bar{q}^1:=q^2,\, \bar{q}^2:=q^3,\, \bar{q}^3:=q^4$, the
Birkhoffian has the form
$\bar{\omega}_c=\sum^{3}_{j=1}\bar{Q}_jd\bar{q}^j$, where \ba
\bar{Q}_1(\bar{q},\dot{\bar{q}},\ddot{\bar{q}})&=&
\left[L_1(\dot{\bar{q}}^1+\dot{\bar{q}}^2)+L_3(\dot{\bar{q}^1})\right]\ddot{\bar{q}}^1+
L_1(\dot{\bar{q}}^1+\dot{\bar{q}}^2)\ddot{\bar{q}}^2-
\nonumber\\
&&C_1(f(\bar{q}^1,\bar{q}^2,\bar{q}^3)-\bar{q}^1+\bar{q}^3+const)\nonumber\\
\bar{Q}_2(\bar{q},\dot{\bar{q}},\ddot{\bar{q}})&=&
L_1(\dot{\bar{q}}^1+\dot{\bar{q}}^2)\ddot{\bar{q}}^1+
L_1(\dot{\bar{q}}^1+\dot{\bar{q}}^2)\ddot{\bar{q}}^2+
R_1(\dot{\bar{q}}^2)+\nonumber\\
&&C_2(f(\bar{q}^1,\bar{q}^2,\bar{q}^3)+\bar{q}^2+const)\nonumber\\
\bar{Q}_3(\bar{q},\dot{\bar{q}},\ddot{\bar{q}})&=&L_2(\dot{\bar{q}}^3)\ddot{\bar{q}}^3+C_1(f(\bar{q}^1,\bar{q}^2,\bar{q}^3)-\bar{q}^1+\bar{q}^3+const)
\label{exrlc39}
\ea
 $f:U\subset\mathbf{R}^3\longrightarrow \mathbf{R}^1$ being the unique function such that
$f(\bar{q}_0)=q^1_0$, $q^1_0\in \mathbf{R}$, and
$C_1(f(\bar{q})-\bar{q}^1+\bar{q}^3+const)+C_2(f(\bar{q})+\bar{q}^2+const)+C_3(f(\bar{q}))=0$,
$\forall \bar{q}=(\bar{q}^1,\bar{q}^2,\bar{q}^3)\in U$, with $U$ a
neighborhood of $\bar{q}_0=(\bar{q}_0^1,\bar{q}_0^2,\bar{q}_0^3)$.

We have shown in section 3 that in this case the reduced
Birkhoffian (\ref{exrlc39}) is   \textbf{dissipative} and
\textbf{regular}.\\
The relations (\ref{bir5}) are satisfied for this example, thus,
there exists a function
$\bar{E}_{0_\omega}(\bar{q},\dot{\bar{q}})$ such that
(\ref{dissip'}) is fulfilled, with the dissipative 1-form given by
\be D=R_1(\dot{\bar{q}}^2)d\bar{q}^2\label{dissip3} \ee We
calculate  \be \textrm{det}\left[\f{\pa\bar{Q}_j}{\pa
\ddot{\bar{q}}^i}\right]_{i,j=1,2,3}= \left|
\begin{array}{ccc}
L_1(\dot{\bar{q}}^1+\dot{\bar{q}}^2)+L_3(\dot{\bar{q}^1})&L_1(\dot{\bar{q}}^1+\dot{\bar{q}}^2)&0\\
L_1(\dot{\bar{q}}^1+\dot{\bar{q}}^2)&L_1(\dot{\bar{q}}^1+\dot{\bar{q}}^2)&0\\
0&0&L_2(\dot{\bar{q}}^3)
\end{array}
\right| \ee Because $L_1,L_2, L_3:\mathbf{R}\longrightarrow
\mathbf{R}\backslash \{0\}$, the determinant above is different
from zero, then, the reduced Birkhoffian given by (\ref{exrlc39})
is regular. $\quad\blacksquare$

 If the nonlinear resistor is voltage controlled, that
is, \be \textsc{i}_{[1]}=\mathfrak{R}_1(v_1) \ee
 we obtain, instead  of
(\ref{exrlc36}), the following implicit Birkhoffian \ba
\hspace{-1cm}\left\{
\begin{array}{ll}
Q_1(q,\dot{q},\ddot{q})=C_1(q^1-q^2+q^4+const)+C_2(q^1+q^3+const)+
C_3(q^1)\\
Q_2(q,\dot{q},\ddot{q})=\left[L_1(\dot{q}^2+\dot{q}^3)+
L_3(\dot{q}^2)\right]\ddot{q}^2+L_1(\dot{q}^2+\dot{q}^3)\ddot{q}^3-\\
\hspace{2.5cm}C_1(q^1-q^2+q^4+const)\\
Q_3(q,\dot{q},\ddot{q})=L_1(\dot{q}^2+\dot{q}^3)\ddot{q}^2+L_1(\dot{q}^2+\dot{q}^3)\ddot{q}^3+C_2(q^1+q^3+const)+v_1\\
Q_4(q,\dot{q},\ddot{q})=L_2(\dot{q}^4)\ddot{q}^4+C_1(q^1-q^2+q^4+const)\\
\\
\dot{q}^3=\mathfrak{R}_1(v_1)
\end{array}\right.
 \label{exrlc36'} \ea
We suppose that the nonlinear voltage controlled resistor
satisfies  $\mathfrak{R}_1(x)x>0$, for any $x\neq 0$. Thus, as in
 the case of current controlled resistors, there exists a
function $E_{0_\omega}(q,\dot{q})$ given by (\ref{exrlc37}), such
that (\ref{dissip}) is satisfied with the dissipative implicit
1-form \be\left\{
\begin{array}{ll}
D(q, \, \dot{q})=
v_1dq^3,\\
\\
\dot{q}^3=\mathfrak{R}_1(v_1)
\end{array}
\right.  \ee Therefore, the Birkhoffian (\ref{exrlc36'})
 is \textbf{dissipative}. The Birkhoffian
(\ref{exrlc36'}) is  also \textbf{\textit{not} regular}.
$\quad\blacksquare$

Example 2)  This example is based  on the following oriented
connected graph

\vspace{0.5cm}
\begin{center}
\scalebox{0.40}{\includegraphics{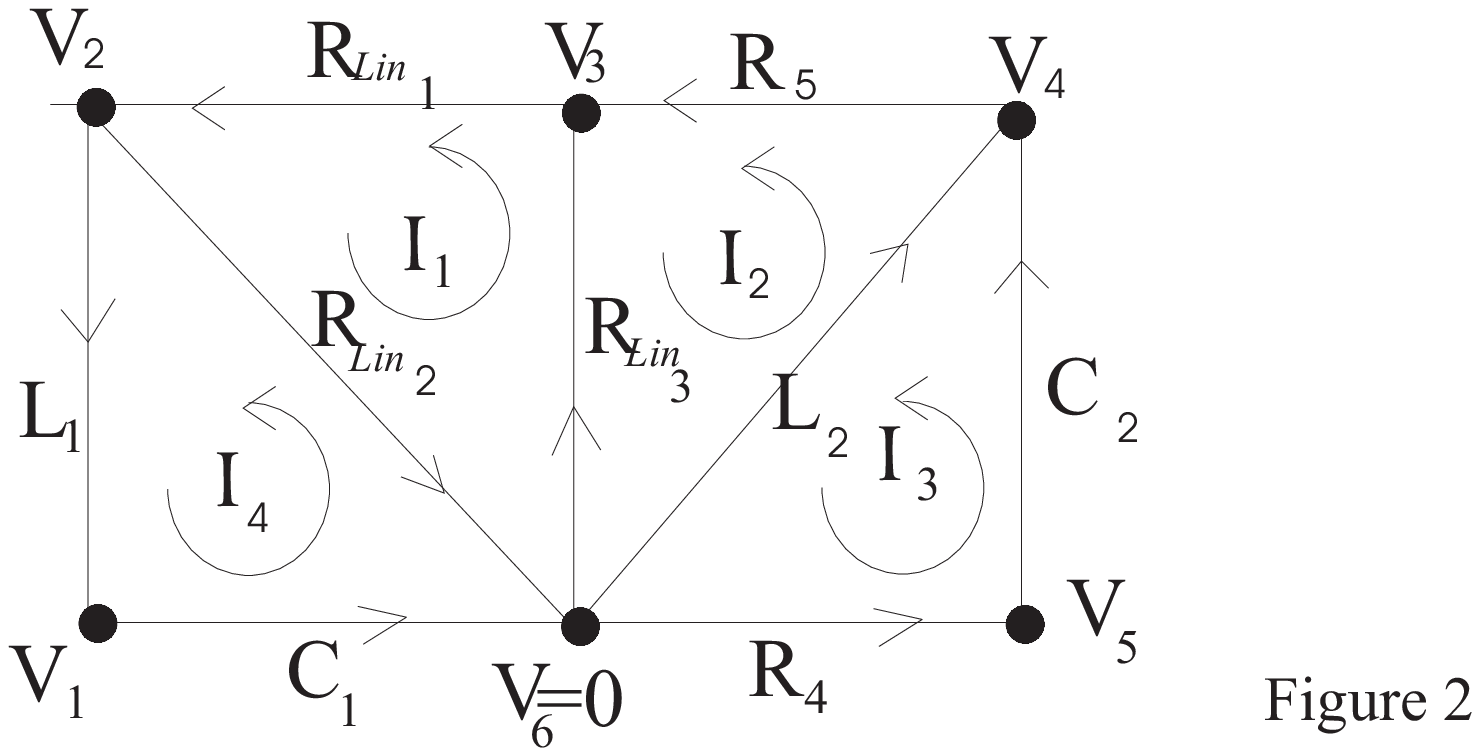}}
\end{center}

\vspace{0.5cm} We have $r=5$, $k=2$, $p=2$, $n=5$, $m=4$ , $b=9$.
We choose the reference node to be $V_6$ and the current
directions as indicated in Figure 2. We cover the associated graph
with the loops $I_1, \, I_2,\, I_3,\, I_4$. Let
$V=(V_1,V_2,V_3,V_4,V_5)\in \mathbf{R}^5$ be the vector of node
voltage values, $\textsc{i}=(\textsc{i}_{[\Gamma]},
\textsc{i}_{(a)},\textsc{i}_\alpha)\in \mathbf{R}^5\times
\mathbf{R}^2\times \mathbf{R}^2$ be
 the vector of branch current  values and
 $v=(v_{[\Gamma]},v_({a}),v_{\alpha})\in \mathbf{R}^5\times \mathbf{R}^2\times\mathbf{R}^2$
 be the vector of branch  voltage
 values.\\
The branches in Figure 2 are  labelled as follows: the first, the
second and the third branch are the linear resistive branches,
 $\textsc{r}_{lin_1}$,  $\textsc{r}_{lin_2}$,  $\textsc{r}_{lin_3}$,
 the fourth and the fifth branch are the nonlinear resistive branches
 $\textsc{r}_{4}$,  $\textsc{r}_{5}$,
 the second, the
sixth and the seventh branch are the inductive branches
$\textsc{l}_1$, $\textsc{l}_2$,  and the last two branches are the
capacitor branches $\textsc{c}_1$, $\textsc{c}_2$. The incidence
and loop matrices, $B\in \mathfrak{M}_{95}(\mathbf{R})$ and
 $A\in \mathfrak{M}_{94}(\mathbf{R})$, write as
\be B=\left(
\begin{array}{ccccc}
0&-1& 1& 0& 0 \\
0& 1& 0& 0& 0 \\
0& 0& -1& 0& 0\\
0& 0& 0&  0&-1\\
0& 0&-1& 1& 0\\
-1& 1& 0& 0& 0\\
 0& 0& 0& 1&0\\
1& 0& 0& 0& 0\\
 0& 0& 0& -1&1
\end{array}
\right),\quad  \quad A=\left(
\begin{array}{cccc}
 1& 0 & 0 & 0 \\
 1& 0 & 0 & -1\\
 1& -1&  0& 0\\
 0& 0 & 1 & 0\\
 0& 1 & 0 & 0 \\
 0& 0 & 0 & 1 \\
 0& -1 & 1 & 0\\
 0& 0& 0&  1\\
 0& 0& 1& 0
\end{array}
\right)
 \label{exrlc21}
\ee One has rank$(B)=5$,  rank$(A)=4$.\\
Except from the resistors in the first loop which are considered
linear, all  devices are nonlinear and are
 described by the relation (\ref{2}), (\ref{4}), (\ref{2rezistor}). We suppose
that $\textrm {\scriptsize R}_1,\, \textrm {\scriptsize R}_2,\,
\textrm {\scriptsize R}_3 > 0$,
 are distinct constants, $C_1,C_2:\mathbf{R}\rightarrow \mathbf{R}\backslash \{0\}$,
 $L_1,L_2:\mathbf{R}\rightarrow \mathbf{R}\backslash \{0\}$, smooth invertible functions and
$R_4,R_5:\mathbf{R}\rightarrow \mathbf{R}$ smooth functions such
that, for any $x\neq 0$ \be R_4(x)x>0, \quad R_5(x)x>0 \label{R45}
\ee The first set of equations (\ref{7})  has the form \be \left\{
\begin{array}{llllllll}
-\textsc{i}_{(1)}+\dot{\textsc{q}}_1=0\\
-\textsc{i}_{[1]}+\textsc{i}_{[2]}+\textsc{i}_{(1)}=0\\
\textsc{i}_{[1]}-\textsc{i}_{[3]}-
\textsc{i}_{[5]}=0\\
\textsc{i}_{[5]}+\textsc{i}_{(2)}-\dot{\textsc{q}}_2=0\\
-\textsc{i}_{[4]}+\dot{\textsc{q}}_2=0
\end{array}
\right.\label{exrlc22} \ee
 The relations (\ref{var}), (\ref{cs}) read as follows
for this example \be
\textsc{i}_{[\Gamma]}:=\dot{\textsc{q}}_{[\Gamma]}, \quad
\Gamma=1,...,5,\quad \textsc{i}_{(a)}:=\dot{\textsc{q}}_{(a)},
\quad a=1,2\label{varrlc2} \ee \be x^1:=\textsc{q}_{[1]},..., \,
x^{5}:=\textsc{q}_{[5]},\,
x^{6}:=\textsc{q}_{(1)},\,x^{7}:=\textsc{q}_{(2)},\,
x^{8}:=\textsc{q}_{1},\, x^{9}:=\textsc{q}_{2}\label{csrlc2} \ee
Using the equations from (\ref{exrlc22}), we define the
4-dimensional affine-linear configuration  space $M_c$. In view of
the notations (\ref{varrlc2}), (\ref{csrlc2}), we integrate these
5 equations and solving them  in terms of $4$ variables, we
obtain, for example, $ x^2=x^{1}-x^{8}+const,\,
x^3=x^1-x^{5}+const,\, x^4=x^9+const,\, x^6=x^8+const,\,
x^7=-x^5+x^9+const.$ Thus, a coordinate system on $M_c$ is given
by \be q^1:=x^1, q^2:=x^5, q^3:=x^8,q^4:=x^9 \label{qcoord} \ee
The matrix of constants $\mathcal{N}$  from (\ref{10'}) is given
by \\
$ \mathcal{N}=\left(
\begin{array}{cccc}
 1& 0& 0& 0\\
 1& 0&-1& 0\\
 1&-1& 0& 0\\
 0& 0& 0& 1\\
 0& 1& 0& 0\\
 0& 0& 1& 0\\
 0&-1& 0& 1\\
0&0&1&0\\
0&0&0&1
\end{array}
\right)$.
 Therefore, in terms of the $q$-coordinates
(\ref{qcoord}), we define the Birkhoffian
$\omega_c=Q_1(q,\dot{q},\ddot{q})dq^1+Q_2(q,\dot{q},\ddot{q})dq^2+Q_3(q,\dot{q},\ddot{q})dq^3+Q_4(q,\dot{q},\ddot{q})dq^4$
of $M_c$ as in (\ref{bir}), with (\ref{bir2}), (\ref{bir3}),
(\ref{bir1}), that is,  \ba
&&\hspace{-0.7cm}Q_1(q,\dot{q},\ddot{q})=(\textrm {\scriptsize
R}_1+\textrm {\scriptsize R}_2+\textrm {\scriptsize
R}_3)\dot{q}^1-
\textrm {\scriptsize R}_3\dot{q}^2-\textrm {\scriptsize R}_2\dot{q}^3\textrm\nonumber\\
&&\hspace{-0.7cm}Q_2(q,\dot{q},\ddot{q})=
L_2(-\dot{q}^2+\dot{q}^4)\ddot{q}^2-L_2(-\dot{q}^2+\dot{q}^4)\ddot{q}^4-\textrm
{\scriptsize R}_3\dot{q}^1+ \textrm {\scriptsize R}_3\dot{q}^2
+R_5(\dot{q}^2)\nonumber\\
&&\hspace{-0.7cm}Q_3(q,\dot{q},\ddot{q})=
L_1(\dot{q}^3)\ddot{q}^3-\textrm {\scriptsize R}_2\dot{q}^1+
\textrm {\scriptsize R}_2\dot{q}^3
+C_1(q^3)\nonumber\\
&&\hspace{-0.7cm}Q_4(q,\dot{q},\ddot{q})=-L_2(-\dot{q}^2+\dot{q}^4)\ddot{q}^2+L_2(-\dot{q}^2+\dot{q}^4)\ddot{q}^4+
R_4(\dot{q}^4)+ C_2(q^4) \label{exrlc26} \ea The Birkhoffian
(\ref{exrlc26}) is \textbf{dissipative} and \textbf{\textit{not}
regular}.\\
 Indeed, there exists a smooth function
$E_{0_{\omega}}:TM\to \mathbf{R}$ of the form (\ref{energie}),
that is, \ba \hspace{-0.7cm}E_{0_\omega}(q,\dot{q})&=& \int
L_1(\dot{q}^3)\dot{q}^3d\dot{q}^3+ \int
\widetilde{L}_2(\dot{q}^2,\dot{q}^4)(-\dot{q}^2+\dot{q}^4)(-d\dot{q}^2+d\dot{q}^4)-
\nonumber\\
&& \int \int
\widetilde{L}_2'(\dot{q}^2,\dot{q}^4)(-\dot{q}^2+\dot{q}^4)d\dot{q}^2d\dot{q}^4+
\int \int \widetilde{L}_2(\dot{q}^2,\dot{q}^4)d\dot{q}^2d\dot{q}^4+\nonumber\\
&&\int C_1(q^3)dq^3+\int C_2(q^4)dq^4 \label{exrlc27} \ea
 such that (\ref{dissip}) is satisfied with the dissipative 1-form defined
 by
\ba \hspace{-1cm}D&=&\left[\dot{q}^1\textrm {\scriptsize R}_1+
(\dot{q}^1-\dot{q}^3)\textrm {\scriptsize R}_2+
(\dot{q}^1-\dot{q}^2)\textrm {\scriptsize R}_3
\right]dq^1+\left[-(\dot{q}^1-\dot{q}^2)\textrm {\scriptsize
R}_3+\right.\nonumber\\
 && \left. R_5(\dot{q}^2)\right] dq^2-(\dot{q}^1-\dot{q}^3)\textrm
{\scriptsize R}_2dq^3+R_4(\dot{q}^4)dq^4\label{D} \ea In view of
the assumptions (\ref{R45})  and of $\textrm {\scriptsize
R}_1,\textrm {\scriptsize R}_2,\textrm {\scriptsize R}_3>0$, we
get \be \hspace{-0.7cm}D_j(q,\dot{q})\dot{q}^j=\textrm
{\scriptsize R}_1(\dot{q}^1)^2+\textrm {\scriptsize
R}_2(\dot{q}^1-\dot{q}^3)^2+\textrm {\scriptsize
R}_3(\dot{q}^1-\dot{q}^2)^2+
R_4(\dot{q}^4)\dot{q}^4+R_5(\dot{q}^2)\dot{q}^2 >0 \ee Therefore,
the vertical 1-form in (\ref{D}) is  dissipative.\\
We are in the case where the  considered circuit has one loop
which contains only resistors. We note that for the Birkhoffian
(\ref{exrlc26}), the first row  of the matrix $\left[\f{\pa
Q_j}{\pa \ddot{q}^i}\right]_{i,j=1,...,4}$ contains only  zeros.
Therefore, $\textrm{det}\left[\f{\pa Q_j}{\pa
\ddot{q}^i}\right]_{i,j=1,...,4}=0$ and the Birkhoffian
(\ref{exrlc26}) is not regular.

\noindent Using the first relation in (\ref{exrlc26}), we now
define $\hat{M}_{c}\subset M_c$ by \be\hspace{-0.7cm}
\hat{M}_{c}=\{q=(q^1,q^2,q^3,q^4)\in M_c/\, (\textrm {\scriptsize
R}_1+\textrm {\scriptsize R}_2+\textrm {\scriptsize
R}_3)q^1-\textrm {\scriptsize R}_3q^2-\textrm {\scriptsize
R}_2q^3+c_1=0 \}\label{reduced} \ee where $c_1$ is a real constant. \\
On the reduced configuration space $\hat{M}_c$, in the coordinate
system given by  $\hat{q}^1:=q^2,\, \hat{q}^2:=q^3,\,
\hat{q}^3:=q^4$, the Birkhoffian has the form
$\hat{\omega}_c=\hat{Q}_1d\hat{q}^1+\hat{Q}_2d\hat{q}^2+\hat{Q}_3d\hat{q}^3$
 \ba
\hspace{-0.7cm}\hat{Q}_1(\hat{q},\dot{\hat{q}},\ddot{\hat{q}})&=&
L_2(-\dot{\hat{q}}^1+\dot{\hat{q}}^3)\ddot{\hat{q}}^1-L_2(-\dot{\hat{q}}^1+\dot{\hat{q}}^3)\ddot{\hat{q}}^3+\left(\mathfrak{C}_2+\mathfrak{C}_3\right)\dot{\hat{q}}^1-
\mathfrak{C}_3\dot{\hat{q}}^2 +R_5(\dot{\hat{q}}^1)
\nonumber\\
\hspace{-0.7cm}\hat{Q}_2(\hat{q},\dot{\hat{q}},\ddot{\hat{q}})&=&
L_1(\dot{\hat{q}}^2)\ddot{\hat{q}}^2-\mathfrak{C}_3\dot{\hat{q}}^1+
\left(\mathfrak{C}_1+\mathfrak{C}_3\right)\dot{\hat{q}}^2
+C_1(\hat{q}^2)
\nonumber\\
\hspace{-0.7cm}\hat{Q}_3(\hat{q},\dot{\hat{q}},\ddot{\hat{q}})&=&
-L_2(-\dot{\hat{q}}^1+\dot{\hat{q}}^3)\ddot{\hat{q}}^1+L_2(-\dot{\hat{q}}^1+\dot{\hat{q}}^3)\ddot{\hat{q}}^3+
R_4(\dot{\hat{q}}^3)+ C_2(\hat{q}^3) \label{exrlc29} \ea where we
denote the constants $\mathfrak{C}_1:=\f {\textrm {\scriptsize
R}_1\textrm {\scriptsize R}_2}{\textrm {\scriptsize R}_1+\textrm
{\scriptsize R}_2+\textrm {\scriptsize R}_3}$,
$\mathfrak{C}_2:=\f{\textrm {\scriptsize R}_1\textrm {\scriptsize
R}_3}{\textrm {\scriptsize R}_1+\textrm {\scriptsize R}_2+\textrm
{\scriptsize R}_3}$, $\mathfrak{C}_3:=\f{\textrm {\scriptsize
R}_2\textrm {\scriptsize R}_3}{\textrm {\scriptsize R}_1+\textrm
{\scriptsize R}_2+\textrm {\scriptsize R}_3}$.

\noindent As we have stated in  section 3, the Birkhoffian given
by (\ref{exrlc29}) is still \textbf{dissipative}. The function
$\hat{E}_{0_\omega}(\hat{q},\dot{\hat{q}})$ has the same form
(\ref{exrlc27}) written in the coordinates $\hat{q}$. The relation
(\ref{dissip1}) is satisfied with the dissipative 1-form defined
by \ba
\hspace{-0.7cm}\hat{D}&=&\left[\dot{\hat{q}}^1\mathfrak{C}_2+(\dot{\hat{q}}^1-\dot{\hat{q}}^2)\mathfrak{C}_3+
R_5(\dot{\hat{q}}^1)\right]d\hat{q}^1+\left[\dot{\hat{q}}^2\mathfrak{C}_1-(\dot{\hat{q}}^1-\dot{\hat{q}}^2)\mathfrak{C}_3\right]d\hat{q}^2+
\nonumber\\
&&R_4(\dot{\hat{q}}^3)d\hat{q}^3 \ea The vertical 1-form above is
 dissipative, as can be seen as follows: For $\textrm
{\scriptsize R}_1>0$, $\textrm {\scriptsize R}_2>0$, $\textrm
{\scriptsize R}_3>0$, we get $\mathfrak{C}_1>0$,
$\mathfrak{C}_2>0$, $\mathfrak{C}_3>0$ and together with
(\ref{R45})  yield \be
\hspace{-1cm}\hat{D}_j(\hat{q},\dot{\hat{q}})\dot{\hat{q}}^j=\mathfrak{C}_2(\dot{\hat{q}}^1)^2+\mathfrak{C}_3(\dot{\hat{q}}^1-\dot{\hat{q}}^2)^2+\mathfrak{C}_1(\dot{\hat{q}}^2)^2+
R_5(\dot{\hat{q}}^1)\dot{\hat{q}}^1+R_4(\dot{\hat{q}}^3)\dot{\hat{q}}^3>0
\ee The Birkhoffian given by (\ref{exrlc29}) is
\textbf{\textit{not} regular}, since  the determinant  \be
\textrm{det}\left[\f{\pa \hat{Q}_j}{\pa
\ddot{\hat{q}}^i}\right]_{i,j=1,2,3}=\left|
\begin{array}{ccc}
\widetilde{L}_2(\dot{\hat{q}})&0&-\widetilde{L}_2(\dot{\hat{q}})\\
0&\widetilde{L}_1(\dot{\hat{q}})&0\\
-\widetilde{L}_2(\dot{\hat{q}})&0&\widetilde{L}_2(\dot{\hat{q}})
\end{array}
\right|= 0 \ee This result does not come unexpected because the
considered network has also a loop which contains only resistors
and capacitors,  formed by $\textsc{r}_4$, $C_2$, $\textsc{r}_5$,
$R_{lin_3}$. In order to regularize the Birkhofiian
(\ref{exrlc29}), we introduce an inductor in series into this
loop, described by the following relation between the current and
the voltage:
$v=\mathcal{L}_{1}({\textsc{i}})\f{d{\textsc{i}}}{dt}$,
$\mathcal{L}_{1}:\mathbf{R}\longrightarrow \mathbf{R}\backslash
\{0\}$ being smooth invertible function. This means that this loop
will have one more node and one more branch. The number of
branches for the  graph associated to the circuit, increases by
one, that is, there will be $b=10$ branches, and the number of
nodes increases by one as well, that is, we will have $n=6$. But
the cardinality of a selection of loops which cover the whole
graph remains $m=4$. After  the calculation we arrive at the
reduced configuration space defined by (\ref{reduced}). On the
reduced configuration space $\hat{M}_c$, in the coordinate system
given by $\hat{q}^1:=q^2,\, \hat{q}^2:=q^3,\, \hat{q}^3:=q^4$, the
Birkhoffian  $\hat{\omega}^{ext}_c$ has the components
$\hat{Q}_1(\hat{q},\dot{\hat{q}},\ddot{\hat{q}})$,
$\hat{Q}_2(\hat{q},\dot{\hat{q}},\ddot{\hat{q}})$ given by
(\ref{exrlc29}) and  the expression of
$\hat{Q}_3(\hat{q},\dot{\hat{q}},\ddot{\hat{q}})$ becomes
\ba\hat{Q}_3(\hat{q},\dot{\hat{q}},\ddot{\hat{q}})&=& -
L_2(-\dot{\hat{q}}^1+\dot{\hat{q}}^3)\ddot{\hat{q}}^1+\left[\mathcal{L}_{1}(\dot{\hat{q}}^3)+
L_2(-\dot{\hat{q}}^1+\dot{\hat{q}}^3)\right]\ddot{\hat{q}}^3+\nonumber\\
&& R_4(\dot{\hat{q}}^3)+ C_2(\hat{q}^3) \ea We now calculate \be
\textrm{det}\left[\f{\pa \hat{Q}_j}{\pa
\ddot{\hat{q}}^i}\right]_{i,j=1,2,3}=\left|
\begin{array}{ccc}
\widetilde{L}_2(\dot{\hat{q}})&0&-\widetilde{L}_2(\dot{\hat{q}})\\
0&\widetilde{L}_1(\dot{\hat{q}})&0\\
-\widetilde{L}_2(\dot{\hat{q}})&0&
\widetilde{\mathcal{L}}_{1}(\dot{\hat{q}})+\widetilde{L}_2(\dot{\hat{q}})
\end{array}
\right|\ee Because
$\mathcal{L}_{1},L_1,L_2:\mathbf{R}\longrightarrow
\mathbf{R}\backslash \{0\}$ the determinant above is different
from zero, then, the Birkhoffian  $\hat{\omega}^{ext}_c$ is
regular. $\quad\blacksquare$

\end{document}